\documentclass[a4paper]{amsart}
\usepackage{lmodern,microtype,amsmath,amssymb,amsthm,mathtools,booktabs,enumitem}
\usepackage{tikz}
\usetikzlibrary{automata,decorations,decorations.markings,decorations.pathmorphing,arrows.meta,quotes}
\usepackage[margin=35mm]{geometry}
\usepackage{hyperref}
\usepackage{lineno}
\hypersetup{pdftitle={Automatic partition identities}}

\newtheorem{theorem}{Theorem}[section]
\newtheorem{lemma}[theorem]{Lemma}
\newtheorem{proposition}[theorem]{Proposition}
\newtheorem{corollary}[theorem]{Corollary}
\newtheorem{maintheorem}{Theorem}

\theoremstyle{definition}
\newtheorem{definition}[theorem]{Definition}

\theoremstyle{remark}

\newcommand{\PP}{{\mathcal P}}
\newcommand{\Ups}{\Upsilon}
\newcommand\K{{\mathbb K}}
\newcommand\Q{{\mathbb Q}}
\newcommand\Z{{\mathbb Z}}
\newcommand\A{{\mathbb A}}
\newcommand\I{{\mathbb I}}

\allowdisplaybreaks[1]
\title{Automatic partition identities}
\author{Laurent Bartholdi}
\date{2 November 2026}
\begin{document}
\begin{abstract}
  We encode certain sets of partitions by regular languages of boundary words,
  so that their generating functions become rational series in the
  quantum plane $RU=qUR$. An equality between two partition counts is
  certified by writing the difference of their series as
  $(1-U)a+b(1-R)$, with $a,b$ converging $q$-adically. We show that
  $a,b$ may be chosen as solutions of finite holonomic systems of $q$-difference
  equations computed from the automata, and describe finite
  certificates of their convergence, together with bounded searches for
  such certificates; the searches need not succeed, but a false identity
  is always detected. As an application, we prove the three-row,
  level-two case of Primc's conjecture, which equates counts of
  coloured partitions subject to local multiplicity conditions with
  counts of coloured partitions whose parts lie in prescribed residue
  classes modulo eight.
\end{abstract}
\maketitle

\begin{quote}
  \noindent\textit{on veut seulement être sûr que c'est elle, ne pas se
  tromper sur l'identité [\dots] cet extrait algébrique\dots}
  \hfill--- Albertine disparue
%
\end{quote}

Taking a cue from Marcel Proust,
I propose an algebraic approach to partition identities, via regular grammars. Many
partition identities take the form `$A(n)=B(n)$', where $A(n)$ and $B(n)$
count the partitions of $n$ satisfying constraints given by regular
expressions. I give in Theorem~\ref{thm:A} a method that (if it doesn't run forever) produces algebraically verifiable certificates of such identities. As an application, I prove in Theorem~\ref{p2:thm:main} the three-row, level-two case of
Primc's conjecture.

\section{A teaser: Euler's identity}\label{sec:teaser}
Here is an illustration of how it works. In Euler's identity, the archetypal example, $A(n)$ counts the
partitions of $n$ into odd parts and $B(n)$ those into distinct parts.
Draw the Ferrers diagram of $\lambda\vdash n$ with $\lambda_i$ boxes in
row~$i$, the longest row on top, and read its lower right boundary from
bottom left to top right, writing $R$ for a unit step to the right and
$U$ for a unit step up. Nonempty partitions thus correspond to words
in $R\{R,U\}^*U$; we ignore the empty partition unless stated otherwise. The odd and
distinct constraints correspond to the regular languages
\[
 \mathcal L_{\mathrm o}=R(U+RR)^*U,
 \qquad
 \mathcal L_{\mathrm d}=R(R+UR)^*U,
\]
whose generating functions are the series
\[
 H_{\mathrm o}=R(1-U-R^2)^{-1}U,
 \qquad
 H_{\mathrm d}=R(1-R-UR)^{-1}U
\]
in the skew field of fractions of the quantum plane
$\A=\Q[q]\langle U,R\rangle/(RU-qUR)$. Expanded as sums over words,
these rational expressions define elements of the $q$-adic completion
$\widehat\A$ of $\A$: the outer letters $R$ and $U$ make the sums
converge, see \S\ref{ss:kernels}. This ring is convenient because the
boundary word of $\lambda=(\lambda_1\ge\dots\ge\lambda_\ell>0)\vdash n$
satisfies
\begin{equation}\label{eq:word-weight}
 w(\lambda)=R^{\lambda_\ell}UR^{\lambda_{\ell-1}-\lambda_\ell}U\cdots
 R^{\lambda_1-\lambda_2}U=q^nU^\ell R^{\lambda_1}\qquad\text{in }\A:
\end{equation}
every pair formed by a letter $R$ and a later letter $U$ contributes a
factor $q$, and the $U$ recording a part $j$ has exactly $j$ such
predecessors.

Every $s\in\widehat\A$ has a unique normal form
$\sum_{i,j}s_{i,j}U^iR^j$ with $s_{i,j}\in\Q[[q]]$. Its normal-ordered
evaluation
\[
 \varepsilon(s)=\sum_{i,j}s_{i,j}\in\Q[[q]]
\]
is linear but not multiplicative: $\varepsilon(RU)=q\neq\varepsilon(R)\varepsilon(U)$. The generating functions of odd and of distinct partitions are
$\varepsilon(H_{\mathrm o})$ and $\varepsilon(H_{\mathrm d})$, and Euler's
identity is $\varepsilon(H_{\mathrm o}-H_{\mathrm d})=0$. A
\emph{certificate} of it is a decomposition
\[
 H_{\mathrm o}-H_{\mathrm d}=(1-U)a+b(1-R),\qquad a,b\in\widehat\A.
\]
Left multiplication by $U$ and right multiplication by $R$ shift normal
exponents without changing coefficients, so
$\varepsilon((1-U)a)=\varepsilon(b(1-R))=0$ by telescoping, and such a
decomposition proves the identity.

The series $a,b$ are typically not rational, but they can be chosen
\emph{holonomic}, that is, defined by finitely many $q$-difference
equations. The \emph{normal symbol} of $s$ is the commutative series
$\sum_{i,j}s_{i,j}x^iR^j$, and its \emph{marginal} is obtained by
setting $x=1$. Since $U^iR^j\cdot U=q^jU^{i+1}R^j$, right
multiplication by $U$ acts on symbols by $g\mapsto x\,\Ups g$, where
\[
 (\Ups g)(x,R)=g(x,qR),
\]
and right multiplication by $R$ acts by $g\mapsto Rg$; on marginals,
$U$ acts by $\Ups$ alone. Since $\Ups R=qR\Ups$, the assignment
$U\mapsto\Ups$, $R\mapsto R$ embeds $\A$ into the opposite of the
algebra of $q$-difference operators in $R$. Consequently, a series
defined by a linear system over $\A$, such as $H_{\mathrm o}$ or
$H_{\mathrm d}$, has a symbol defined by a linear system of
$q$-difference equations in $R$, with $x$ as a parameter. A certificate
then consists of such systems for $a$ and $b$, finitely many values
selecting their solutions, and finite data guaranteeing that their
coefficients tend $q$-adically to $0$. For instance, write $\nu$ for
the $q$-adic valuation, and call a series \emph{$q$-integral} if it has
no negative powers of $q$. If $b=q\Ups b'$ with $b'$ $q$-integral,
that is, $b'\in\Q[[q,R]]$, then $[R^n]b=q^{n+1}[R^n]b'$ has valuation
at least $n+1$, so $b\in\widehat\A$. A holonomic kernel for $b'$,
together with finite data certifying its integrality, is such a
guarantee.

For Euler's identity, moving the outer letters inwards gives
$H_{\mathrm o}=qU(1-qU-q^2R^2)^{-1}R$ and
$H_{\mathrm d}=qU(1-qR-q^2UR)^{-1}R$. Hence
$H_{\mathrm o}-H_{\mathrm d}$ has normal symbol
\[
 D(x,R)=qR(m_{\mathrm o,x}-m_{\mathrm d,x}),\qquad
 \begin{cases}
 (1-qx\Ups-q^2R^2)m_{\mathrm o,x}=x,\\
 (1-qR-q^2xR\Ups)m_{\mathrm d,x}=x.
 \end{cases}
\]
Both operators are congruent to $1$ modulo $q$, so these equations have
unique solutions in $\Q[x,R][[q]]$. The canonical $a,b$ are given by
\[
 (1-x)a=D(x,R)-D(1,R),\qquad (1-R)b=D(1,R),
\]
and are holonomic. It remains to certify their convergence in
$\widehat\A$. The integral kernels (\S\ref{ss:conv-cert})
\[
 u=1+R^2\alpha(\Ups)u,\quad v=1+R^2\beta(\Ups)v,\qquad
 \alpha(U)=(1-qU)^{-1},\quad\beta(U)=(1+qU)(1+q^2U),
\]
select $u,v\in\Q[[q,R]]$, and $\chi(U)=1-q^2U$ is a \emph{dilation witness}
for $\theta=\beta/\alpha$:
\[
 \theta(U)\chi(U)=(1-q^2U^2)(1-q^4U^2)=\chi(U^2)\chi(q^2U^2).
\]
Thus $u-v=q^2\Ups w$ with
$w\in\Q[[q,R]]$, see Proposition~\ref{prop:dilation-checker}, and the equations of $u$, $v$, $m_{\mathrm o,1}$ and $m_{\mathrm d,1}$ give
\begin{equation}\label{eq:euler-direct}
 b=\bigl(1+R(1-q\Ups)^{-1}\bigr)u-\bigl(1+R(1+q\Ups)\bigr)v=q\Ups b'
\end{equation}
with $b'=qw+R(1-q\Ups)^{-1}(w+\Ups v)\in\Q[[q,R]]$. Thus
$b\in\widehat\A$; moreover $a\in\widehat\A$ by
Theorem~\ref{thm:A}(2), and the certificate is complete.
Corollary~\ref{cor:euler-integral} gives the details, and
\S\ref{ss:glaisher} applies the same mechanism to Glaisher's identity,
of which Euler's is the case $\lambda=2$.

\section{Main results}
The scheme of \S\ref{sec:teaser} applies to any two partition languages
in $R\{R,U\}^*U$ recognized by unambiguous finite automata, and more
generally to weighted automata recording colours or other statistics; a
boundary word is then counted with the total weight of its accepting
paths. Section~\ref{ss:compile} explains how classical constraints are
compiled into such automata. Stripping the first $R$ and the last $U$
yields interior automata with $s_1$ and $s_2$ states. They give a
\emph{normalized realization}
\[
 w\mathsf Q=r,\qquad S=we,\qquad\mathsf Q\equiv \I\pmod q,
\]
of dimension $s=s_1+s_2$, of the difference $S$ of the two series; see
\S\ref{ss:kernels}. We write $\K=\Q(q)$, and use the Pochhammer symbols
$(a;q)_n=\prod_{j=0}^{n-1}(1-aq^j)$ and
$(a_1,\ldots,a_t;q)_n=\prod_i(a_i;q)_n$, also for $n=\infty$ whenever
the product converges $q$-adically.

\begin{maintheorem}\label{thm:A}
Let $H_1,H_2$ be the path series \textup{(\S\ref{ss:languages})} of
finite weighted automata with weights in $\Q(q)$ regular at $q=0$, all of whose
accepted words lie in $R\{R,U\}^*U$. Let $S=H_1-H_2$ have normal symbol
$D(x,R)$ and a normalized realization of dimension $s$, and put
\[
 a(x,R)=\frac{D(x,R)-D(1,R)}{1-x},\qquad
 b(R)=\frac{D(1,R)}{1-R}=\sum_{n\ge0}b_nR^n.
\]
\begin{enumerate}
\item The series $a$ and $b$ are components of the unique
  solutions of explicit square inhomogeneous systems of linear
  $q$-difference equations, of sizes $2s+1$ and $s+1$, computed from the
  automata. In particular, $a$ is holonomic in
  $(x,R)$, $b$ is holonomic in $R$, and every coefficient $b_n$ lies in $\K$.
\item The coefficients of $a$ satisfy $\nu(a_{i,j})\ge i+j$; in
  particular $a\in\widehat\A$.
\item The following are equivalent:
  \begin{enumerate}
  \item $\varepsilon(H_1)=\varepsilon(H_2)$;
  \item $b_n\to0$ $q$-adically;
  \item $\nu(b_n)\ge n+1$ for all $n$;
  \item $b=q\Ups\phi$ for some $\phi\in\Q[[q,R]]$;
  \item $b\in\Q[R][[q]]$, so that $S=(1-U)a+b(1-R)$ in $\widehat\A$.
  \end{enumerate}
  In~(d), $\phi=q^{-1}b(R/q)$, and $\phi$ is holonomic. In~(e),
  $(a,b)$ is the only pair with $S=(1-U)a+b(1-R)$, $a\in\widehat\A$ and
  $b\in\Q[R][[q]]$.
\end{enumerate}
\end{maintheorem}

Theorem~\ref{thm:A} is proved in~\S\ref{ss:proof-A}. It reduces a
partition identity to the integrality of a single holonomic series
$\phi$, defined by finite data. The formal boundaries $a,b$ exist
whether or not the identity holds, and a false identity is detected
after computing finitely many $b_n$: since
$\nu(\varepsilon(S)-b_n)\ge n+1$ unconditionally
(Lemma~\ref{lem:area}), the identity is false if and only if
$\nu(b_n)<n+1$ for some $n$. For a true identity, it remains to prove integrality of $\phi$. A \emph{direct certificate} (\S\ref{ss:conv-cert}) does
so by finitely many exact checks on holonomic series; alternatively, any
independent proof of the identity yields~(3d)
(Corollary~\ref{cor:indirect}).

This separates discovery from verification, in the spirit of
Wilf--Zeilberger certificates~\cite{WZ90,WZ92,Zeilberger}. Certificates
are found by finite linear-algebra searches (\S\ref{ss:searches} and
Appendix~\ref{app:searches}) for integral kernels, dilation witnesses,
finite recurrences, telescoping certificates and quadratic invariants,
and then checked by exact algebraic manipulations. Since the quantum plane records
area, many combinatorial descriptions --- allowed residues, controlled
differences between coloured parts, forbidden patterns --- compile into
rational series to which Theorem~\ref{thm:A} applies
(\S\ref{ss:compile}).

Many partition identities equate restrictions on differences or
multiplicities with restrictions on the residue classes of the parts.
Classical examples include Euler's identity between distinct and odd
parts (modulo~$2$)~\cite[\S326]{Euler1748}, the Rogers--Ramanujan
identities (modulo~$5$)~\cite{Rogers1894}, Schur's theorem
(modulo~$6$)~\cite{Schur1926}, and the G\"ollnitz--Gordon identities
(modulo~$8$)~\cite{Gollnitz1967,Gordon1965}. Later examples include
Capparelli's identity (modulo~$12$), conjectured in his 1988
thesis~\cite{Capparelli1988} and first proved by
Andrews~\cite{Andrews1994}, and Nandi's three identities
(modulo~$14$), conjectured in 2014~\cite[\S8.1]{Nandi2014} and proved by
Takigiku and Tsuchioka~\cite{TakigikuTsuchioka2024}. These examples
suggest regular languages as a common framework; Nandi's difference
conditions already involve forbidden patterns of unbounded length.
Primc's coloured partition conjectures~\cite{Primc} belong to the same
tradition; our main application, Theorem~\ref{p2:thm:main} below, is
their three-row case at level two, with modulus~$8$.

\subsection*{Primc's identities at level two}\label{ss:primc}
We prove the three-row, level-two case of Primc's
conjecture~\cite[Conjecture~2.1]{Primc}. It equates two classes of
coloured partitions: one defined by restrictions on neighbouring
multiplicities, the other by congruence classes modulo eight.

For the first class, every positive even part has two colours, called
\emph{top} and \emph{bottom}, and every odd part has one colour, called
\emph{middle}. Write $\tau_n,\beta_n$ for the multiplicities of the top
and bottom parts $2n$, and $\mu_n$ for the multiplicity of the part
$2n-1$. Fix nonnegative integers $k_0,k_1,k_2$ of sum $2$, the
\emph{level}, and set $\tau_0=k_2$, $\beta_0=k_0$. The partition is
\emph{admissible} if
\begin{equation}\label{p2:eq:local}
 \max(\tau_{n-1},\tau_n)+\mu_n+
 \max(\beta_{n-1},\beta_n)\le2\qquad(n\ge1).
\end{equation}
Equivalently, each of the four sums obtained by choosing one of
$\tau_{n-1},\tau_n$, one of $\beta_{n-1},\beta_n$, and $\mu_n$ is at most
two; in Primc's description by three-row arrays, these are the sums
along the four downward paths through $\mu_n$
(see~\S\ref{ss:compile}). The value
$k_1=2-k_0-k_2$ imposes no further condition. There are six boundary
conditions; write $\Pi_{k_0,k_1,k_2}(q)$ for the generating function of
the corresponding admissible partitions, including the empty one.

\begin{maintheorem}\label{p2:thm:main}
For each $N\ge0$, the number of admissible coloured partitions of $N$
with the boundary condition in the first column equals the number of
coloured partitions of $N$ with the following allowed residues modulo
eight:
\[
\begin{array}{c|c|c}
 (k_0,k_1,k_2)&\text{residues with one colour}&
                \text{residues with two colours}\\\hline
 (2,0,0)&2,3,5,6&4\\
 (1,1,0)&1,2,6,7&4\\
 (1,0,1)&3,5&2,6\\
 (0,2,0)&1,7&2,6.
\end{array}
\]
Each allowed colour has unrestricted multiplicity, and all other
residues are forbidden. The cases $(0,0,2)$ and $(0,1,1)$ follow by
interchanging the top and bottom colours.
\end{maintheorem}

For example, for the boundary $(2,0,0)$ the local
rule~\eqref{p2:eq:local} forbids the part $1$ and the top-coloured part
$2$, and the theorem states
\[
 \Pi_{2,0,0}(q)
   =\frac1{(q^2,q^3,q^4,q^4,q^5,q^6;q^8)_\infty}.
\]
Figure~\ref{fig:primc-automata} shows the two automata for this
boundary; the tags on the letters $U$ make them unambiguous and are
forgotten in $\A$ (\S\ref{ss:compile}).

Primc's Conjecture~2.1 only asserts that these generating functions are
periodic products; the four products of Theorem~\ref{p2:thm:main} are
those he found numerically~\cite[p.~53]{Primc}. Kanade, Russell,
Tsuchioka and Warnaar proposed an explicit product for every width and
level~\cite[Conjecture~4.2]{KRTW}, whose case $n=k=2$, with the same
labelling of the boundary, gives the four products above. They note
that the three-row case appears to be open, and they prove functional
equations~\cite[Propositions~4.4 and~4.5]{KRTW} and quadruple-sum
expressions~\cite[Theorem~4.10]{KRTW} for the four generating
functions, refined by the number of parts, but not their product
evaluations. The level-one case, of modulus~$6$, follows from the
coloured partition theorem~\cite[Theorem~1.2]{JingMisraSavage2001},
and also from a multisum identity of Russell~\cite{Russell} together
with Bressoud's identities; see~\cite[\S4]{KRTW}.

Both partition classes are recognized by small automata: the first
needs only to remember the preceding even multiplicities, and the second
only the part size modulo eight and its colour. The proof
in~\S\ref{sec:primc} reduces the four-state transfer system of the
local conditions to two G\"ollnitz--Gordon systems, by a quadratic
invariant and two rational telescoping certificates, all found by
the searches of~\S\ref{ss:searches}. The classical G\"ollnitz--Gordon
identities and Jacobi's triple product then give the products in the
theorem. This proves $\varepsilon(H_1)=\varepsilon(H_2)$ for the path
series of the two automata, and Corollary~\ref{cor:indirect} turns it
into convergent holonomic certificates for them. The new step is the
reduction of Primc's four-state system to these classical identities.

In~\S\ref{sec:examples} we also treat Glaisher's identity, the two
Rogers--Ramanujan identities, Bessenrodt's refinement of Euler's identity by a global
partition statistic, Berkovich--Uncu identities involving part sizes
and indices, and an identity between two five-state automata,
found by computer search, which fits none of the previous patterns.

\section{Holonomic certificates and their discovery}\label{sec:certificates}

\subsection{Partition languages}\label{ss:languages}
General background on partitions and $q$-series may be found in
Andrews~\cite{Andrews}. In the boundary word $w(\lambda)$
of~\eqref{eq:word-weight}, $R$ raises the current \emph{height}, the number of letters $R$ read so far, and $U$ records a part equal to the height; we write $\ell(\lambda)$ for the number
of parts of $\lambda$ and $|\lambda|$ for their sum.

Regular languages are those described by regular expressions, built
from letters by concatenation, union and Kleene star; equivalently,
those recognized by finite automata, namely finite directed graphs with
an initial state, a set of final states and edges labelled $R$ or $U$,
whose language consists of the labels of paths from the initial state
to a final state. If $\mathsf T$ is the transition matrix of the
automaton, whose entries are sums of letters, and $r,e$ are the
indicator vectors of the initial and final states, then the \emph{path
series} $r(\I-\mathsf T)^{-1}e=\sum_{n\ge0}r\mathsf T^ne$ counts the
accepting paths by their labels; for an unambiguous automaton, it is the
series of its language. Weighted automata allow other coefficients, for
instance multiplicities of colours or markers of statistics.

\subsection{Compiling constraints into automata}\label{ss:compile}
Let $f_h$ be the multiplicity of the part $h$ in $\lambda$. Grouping the
letters of $w(\lambda)$ by height gives
\begin{equation}\label{eq:mult-word}
 w(\lambda)=RU^{f_1}RU^{f_2}\cdots RU^{f_m},\qquad m=\lambda_1,\quad f_m\ge1.
\end{equation}
Reading a boundary word from left to right therefore scans the
multiplicities in increasing order of the parts, and a constraint is
regular as soon as it can be checked during this scan with finite
memory. We make this precise for the constraints occurring in
Rogers--Ramanujan type identities: congruence conditions on the parts,
and difference conditions on the multiplicities, possibly of coloured
parts.

\subsubsection*{Colours and tagged letters}
Following Primc~\cite{Primc}, let
$\mathcal N=\mathcal N_1\cup\dots\cup\mathcal N_r$ be a multiset of
positive integers, the elements of $\mathcal N_c$ having colour $c$. A
coloured partition on $\mathcal N$ is given by multiplicities $f_{h,c}$,
which vanish unless $h\in\mathcal N_c$; we write
$f_h=(f_{h,1},\dots,f_{h,r})$. Encode it by the word
\[
 w(\lambda)=RU_1^{f_{1,1}}\cdots U_r^{f_{1,r}}\,
            RU_1^{f_{2,1}}\cdots U_r^{f_{2,r}}\cdots
            RU_1^{f_{m,1}}\cdots U_r^{f_{m,r}}
\]
over the alphabet $\{R,U_1,\dots,U_r\}$ of \emph{tagged letters}, in
which the letters at each height appear in increasing order of colour.
This encoding is injective. The tags make the automata below
unambiguous; they are forgotten in $\A$, where each $U_c$ becomes $U$,
so that path series count coloured partitions, with multiplicity when
several of them have the same uncoloured boundary word.

\subsubsection*{Congruence conditions}
Suppose that every $\mathcal N_c$ is a union of residue classes modulo
$M$. The encoding above then has regular image: consider the automaton whose states are the pairs $(\rho,j)$, where
$\rho\in\Z/M$ records the current height modulo $M$ and
$j\in\{0,\dots,r\}$ the last colour read at this height ($j=0$ if none),
together with an initial state $\iota$; its final states are the
$(\rho,j)$ with $j\ge1$, and its transitions are
\[
 \iota\xrightarrow{\;R\;}(1,0),\qquad
 (\rho,j)\xrightarrow{\;R\;}(\rho+1,0),\qquad
 (\rho,j)\xrightarrow{\;U_c\;}(\rho,c)\quad
 \text{if $j\le c$ and $\rho\in\mathcal N_c\bmod M$}.
\]
It is deterministic and accepts the encodings of all nonempty coloured
partitions on $\mathcal N$, so its generating function is
\[
 \prod_{c=1}^r\prod_{h\in\mathcal N_c}\frac1{1-q^h}-1
 =\prod_{t=1}^{M}\frac1{(q^t;q^M)_\infty^{c_t}}-1,
\]
where $c_t$ is the number of colours available at the residue of $t$.
Thus every periodic product with nonnegative exponents is, up to the
empty partition, the generating function of a regular language;
Figure~\ref{fig:rr-automata} (right) shows a variant, with one colour
and a separate final state, for $M=5$. A product with negative
exponents, such as Primc's conjectured
$(q^3;q^6)_\infty/(q,q^2,q^4,q^5;q^6)_\infty$ for the boundary
$(0,1,0)$, is first cleared of its numerator: multiplying both sides by
$(q^3;q^6)_\infty^{-1}$ adds, on the side of the difference conditions,
one colour available at the parts $\equiv3\pmod6$ and subject to no
other condition.

\subsubsection*{Difference conditions}
Difference conditions constrain the multiplicities in a window of
consecutive heights. When the multiplicities they involve are bounded,
only finitely many multiplicity vectors occur in a window, and the
conditions define a regular language:

\begin{lemma}[Window conditions]\label{lem:window}
Let $r,M,W\ge1$ and $K\ge0$ be integers. Let
$\mathcal W\subseteq\Z/M\times(\{0,\dots,K\}^r)^W$ contain every tuple
$(\rho,0,\dots,0)$, and fix vectors
$v_{2-W},\dots,v_0\in\{0,\dots,K\}^r$. The coloured partitions such
that
\begin{equation}\label{eq:window}
 (h\bmod M,\,f_{h-W+1},\dots,f_h)\in\mathcal W\qquad\text{for all }h\ge1,
\end{equation}
where $f_h=v_h$ for $h\le0$ and $f_h=0$ for $h>\lambda_1$, are encoded by
a regular language, recognized by a deterministic automaton with at most
$1+M(K+1)^{rW}$ states.
\end{lemma}
\begin{proof}
After reading a prefix ending at height $h$, the automaton records
$h\bmod M$, the vectors $f_{h-W+1},\dots,f_{h-1}$, and the part of $f_h$
read so far. The first letter $R$ leads from the initial state to
height~$1$, with recorded vectors $v_{2-W},\dots,v_0$. A letter $U_c$
increases the $c$-th entry of the current vector; it is rejected if a
larger colour has already been read at this height, or if the entry
would exceed $K$. A letter $R$ is rejected unless the completed window
lies in $\mathcal W$; otherwise the oldest vector is discarded and the
height increases by one. A state reached by a letter $U_c$ is final if
the current window, and the windows at heights $h+1,\dots,h+W-1$
obtained from it by appending zero vectors, lie in $\mathcal W$; all
later windows are zero.
\end{proof}

Only the information needed for later checks has to be stored, and
the automaton can then be minimized. Colours on
which $\mathcal W$ imposes no condition need not be bounded: the
automaton reads their letters without recording them. Conditions of different kinds combine,
since intersections of regular languages are regular. For instance,
Gordon's conditions~\cite{Gordon1961} $f_h+f_{h+1}\le k$ and $f_1\le a$,
with $0\le a\le k$, are~\eqref{eq:window} with $r=M=1$, $W=2$, $K=k$ and the single
initial vector $v_0=k-a$: the initial condition becomes a difference
condition involving an \emph{imaginary multiplicity} at height~$0$. For
$k=1$ and $a=1,0$ these are the Rogers--Ramanujan languages
$\mathcal L_1,\mathcal L_2$ of Figure~\ref{fig:rr-automata} (left).

\subsubsection*{Primc's arrays}
Primc's colours and conditions~\cite{Primc} are of this kind. For
$w=2\ell-1$, the multiset $\mathcal N^{\mathrm{odd}}_w$ has one colour
for each row $i=1,\dots,w$ of an array, row $i$ containing the parts of
the parity of $i+1$: the even parts have $\ell$ colours, and the odd
parts $\ell-1$. Two entries are \emph{adjacent} if they lie in adjacent
rows and their parts differ by one, and a \emph{downward path} is a
sequence of entries, one in each row from top to bottom, any two
consecutive ones being adjacent. The level-$k$ difference conditions
require that the multiplicities along every downward path sum to at
most $k$. The initial conditions add imaginary multiplicities: the
boundary $(k_0,\dots,k_\ell)$, with $k_0+\dots+k_\ell=k$, is the list
$(k_0,k_1,0,k_2,0,\dots,0,k_{\ell-1},k_\ell)$, read from the bottom row
to the top row, placed at the part $0$ in the odd rows and at the part
$-1$ in the even rows.

The parts along a downward path lie in $w$ consecutive integers, and
every multiplicity lies on some downward path. Hence the admissible
coloured partitions are defined by~\eqref{eq:window}, with $r=W=w$,
$M=2$, $K=k$, the set $\mathcal W$ of windows that respect the parities
of the rows and in which every downward path has sum at most $k$, and
initial vectors $v_{2-w},\dots,v_0$ that vanish except for $v_{-1},v_0$,
which carry the imaginary multiplicities. By Lemma~\ref{lem:window},
for every width, level and boundary, the admissible coloured partitions
form a regular language. Every product formula
of~\cite[Conjecture~4.2]{KRTW}, and every product found numerically
in~\cite{Primc}, once cleared of negative exponents as above, therefore
takes the form $\varepsilon(H_1)=\varepsilon(H_2)$ of
Theorem~\ref{thm:A}. For $w=3$, $k\le2$ and for $w=5$, $k=1$, the
resulting automata reproduce the first terms of all products listed
in~\cite[pp.~53--54]{Primc}.

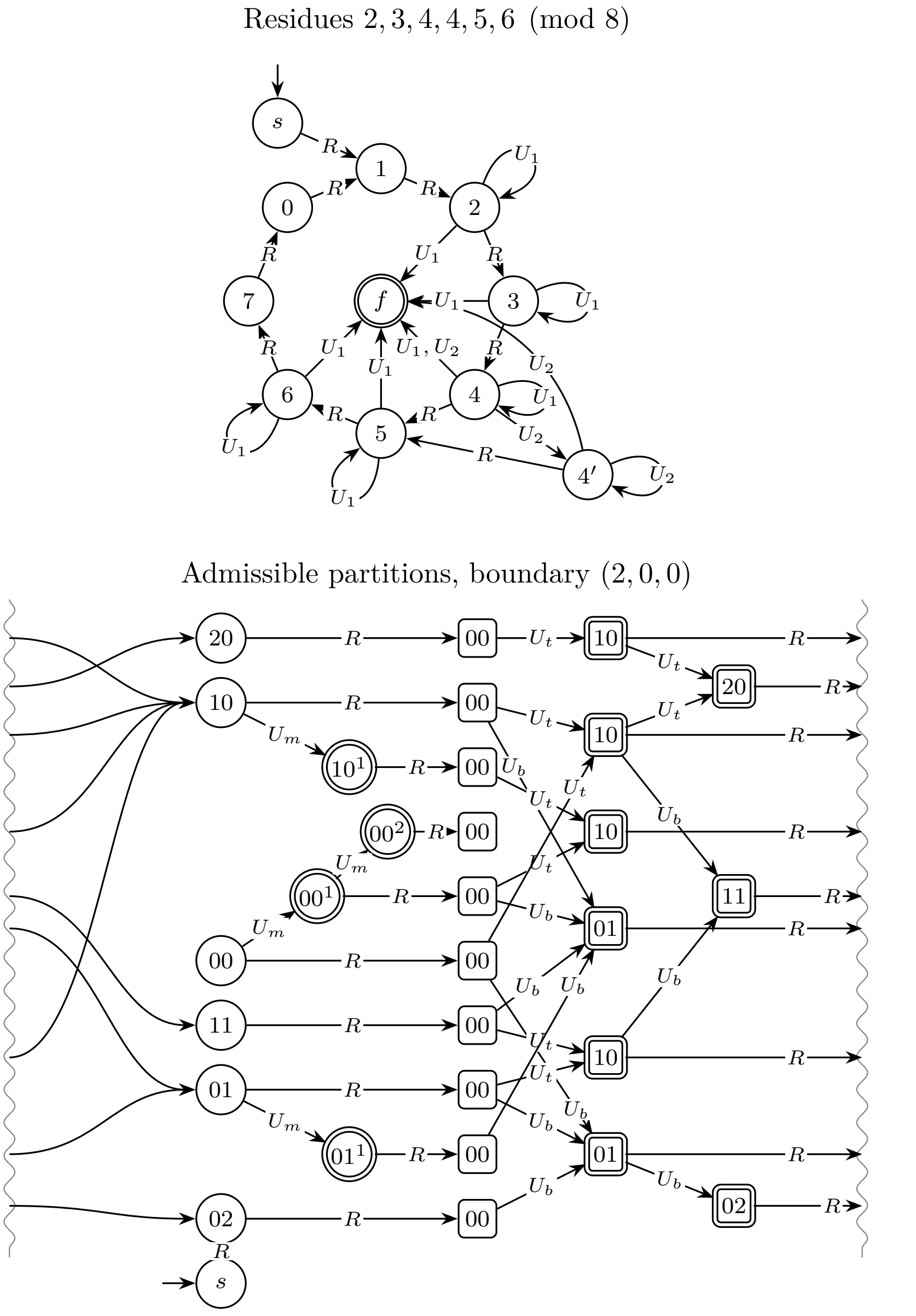
\begin{figure}[!htbp]
\centering
\tikzset{
  pa/.style={
    >=Stealth, semithick,
    initial text={}, initial distance=4mm,
    every state/.style={
      draw,circle,minimum size=17pt,inner sep=0.5pt,font=\footnotesize
    },
    box/.style={
      draw,rectangle,rounded corners=2pt,minimum size=13pt,
      inner sep=1.5pt,font=\footnotesize
    },
    fin/.style={double,double distance=1pt},
    every edge/.style={draw,->},
    every edge quotes/.style={
      auto=false,midway,fill=white,inner sep=0.7pt,font=\scriptsize
    },
    every loop/.style={looseness=6}
  },
  lp/.style 2 args={loop,out=#1,in=#2,min distance=7mm,looseness=9},
  cut/.style={decorate,decoration={snake,amplitude=0.65mm,segment length=4mm},
              thin,gray}
}
Residues $2,3,4,4,5,6\pmod 8$
\par\smallskip
\begin{tikzpicture}[pa]
  \foreach \r/\angle in
    {1/90,2/45,3/0,4/-45,5/-90,6/-135,7/180,0/135}
    \node[state] (r\r) at (\angle:1.6) {$\r$};
  \node[state,initial,initial where=above] (s) at (-1.25,2.15) {$s$};
  \node[state,accepting] (f) at (0,0) {$f$};
  \node[state] (r4p) at (2.5,-2.1) {$4'$};

  \path (s) edge["$R$"] (r1);
  \foreach \j/\k in {1/2,2/3,3/4,4/5,5/6,6/7,7/0,0/1}
    \path (r\j) edge["$R$"] (r\k);
  \path
    (r2) edge[lp={70}{20},"$U_1$"] (r2)
    (r3) edge[lp={25}{-25},"$U_1$"] (r3)
    (r4) edge[lp={20}{-25},"$U_1$"] (r4)
    (r5) edge[lp={-95}{-145},"$U_1$"] (r5)
    (r6) edge[lp={-110}{-160},"$U_1$"] (r6)
    (r4) edge["$U_2$"] (r4p)
    (r4p) edge[lp={25}{-25},"$U_2$"] (r4p)
    (r4p) edge["$R$"] (r5)
    (r4p) edge[bend right=38,"$U_2$" {pos=.38}] (f);
  \foreach \r in {2,3,5,6}
    \path (r\r) edge["$U_1$"] (f);
  \path (r4) edge["{$U_1,U_2$}"] (f);
\end{tikzpicture}

\medskip
Admissible partitions, boundary $(2,0,0)$
\par\smallskip
\begin{tikzpicture}[pa,x=1.55cm,y=0.78cm]
  \draw[cut] (-.65,.6) -- (-.65,-9.6);
  \draw[cut] (6,.6) -- (6,-9.6);
  \node[state] (O20) at (1,0) {$20$};
  \node[state] (O10) at (1,-1) {$10$};
  \node[state] (O00) at (1,-5) {$00$};
  \node[state] (O11) at (1,-6) {$11$};
  \node[state] (O01) at (1,-7) {$01$};
  \node[state] (O02) at (1,-9) {$02$};
  \node[state,initial,initial where=left] (s) at (1,-10) {$s$};
  \node[state,fin] (B10) at (2,-2) {$10^1$};
  \node[state,fin] (B00b) at (2.3,-3) {$00^2$};
  \node[state,fin] (B00a) at (1.75,-4) {$00^1$};
  \node[state,fin] (B01) at (2,-8) {$01^1$};
  \foreach \n/\y in {E20/0,E10/-1,E10b/-2,E00c/-3,E00b/-4,E00/-5,
                     E11/-6,E01/-7,E01b/-8,E02/-9}
    \node[box] (\n) at (3,\y) {$00$};
  \node[box,fin] (Dt22) at (4,0) {$10$};
  \node[box,fin] (Dt8) at (4,-1.5) {$10$};
  \node[box,fin] (Dt11) at (4,-3) {$10$};
  \node[box,fin] (Db6) at (4,-4.5) {$01$};
  \node[box,fin] (Dt21) at (4,-6.5) {$10$};
  \node[box,fin] (Db4) at (4,-8) {$01$};
  \node[box,fin] (Dtt19) at (5,-.75) {$20$};
  \node[box,fin] (Dtb15) at (5,-4) {$11$};
  \node[box,fin] (Dbb23) at (5,-8.8) {$02$};

  \path (s) edge["$R$"] (O02);
  \path
    (O00) edge["$U_m$"] (B00a)
    (B00a) edge["$U_m$"] (B00b)
    (O10) edge["$U_m$"] (B10)
    (O01) edge["$U_m$"] (B01);
  \foreach \a/\b in {O20/E20,O10/E10,B10/E10b,B00b/E00c,B00a/E00b,
                     O00/E00,O11/E11,O01/E01,B01/E01b,O02/E02}
    \path (\a) edge["$R$"] (\b);
  \path
    (E20) edge["$U_t$"] (Dt22)
    (E10) edge["$U_t$"] (Dt8)
    (E10) edge["$U_b$" {pos=.24}] (Db6)
    (E10b) edge["$U_t$"] (Dt11)
    (E00b) edge["$U_t$"] (Dt11)
    (E00b) edge["$U_b$"] (Db6)
    (E00) edge["$U_t$" {pos=.82}] (Dt8)
    (E00) edge["$U_b$" {pos=.84}] (Db4)
    (E11) edge["$U_t$"] (Dt21)
    (E11) edge["$U_b$" {pos=.35}] (Db6)
    (E01) edge["$U_t$"] (Dt21)
    (E01) edge["$U_b$"] (Db4)
    (E01b) edge["$U_b$" {pos=.80}] (Db6)
    (E02) edge["$U_b$"] (Db4)
    (Dt22) edge["$U_t$"] (Dtt19)
    (Dt8) edge["$U_t$"] (Dtt19)
    (Dt8) edge["$U_b$"] (Dtb15)
    (Dt21) edge["$U_b$"] (Dtb15)
    (Db4) edge["$U_b$"] (Dbb23);
  \foreach \a/\y/\b/\target in {
    Dt22/0/O10/-1,Dtt19/-.75/O20/0,Dt8/-1.5/O10/-1,
    Dt11/-3/O10/-1,Dtb15/-4/O11/-6,Db6/-4.5/O01/-7,
    Dt21/-6.5/O10/-1,Db4/-8/O01/-7,Dbb23/-8.8/O02/-9}
  {
    \path (\a) edge["$R$" {pos=.72}] (6,\y);
    \draw[->] (-.65,\y) .. controls (.05,\y) and (.15,\target) .. (\b);
  }
\end{tikzpicture}
\caption{The two automata for boundary $(2,0,0)$.
Top: the residue automaton; $U_1,U_2$ distinguish the colours at~$4$.
Bottom: circles record the preceding even pair $\tau\beta$, with the
number of middle parts as exponent; rectangles record the current
even pair. Double borders mark acceptance. Identify the wavy cuts:
each outgoing $R$-edge resumes at the same vertical position on the
left. The ten plain rectangles also have an $R$-edge to circle $00$
(omitted). Repeated labels have different allowed continuations,
determined by \eqref{p2:eq:local}. All tagged $U$'s become $U$ in
$\A$.}
\label{fig:primc-automata}
\end{figure}

For $w=3$, with the multiplicities $\tau_n,\mu_n,\beta_n$ used in
Theorem~\ref{p2:thm:main}, the array is
\[
\begin{array}{ccccccccccc}
 &\tau_0&&\tau_1&&\tau_2&&\tau_3&&\cdots\\
 \mu_0&&\mu_1&&\mu_2&&\mu_3&&\mu_4&&\cdots\\
 &\beta_0&&\beta_1&&\beta_2&&\beta_3&&\cdots
\end{array}
\]
with columns indexed by the parts $-1,0,1,2,\dots$, and imaginary
multiplicities $\tau_0=k_2$, $\mu_0=k_1$, $\beta_0=k_0$. The four
downward paths through $\mu_n$ give~\eqref{p2:eq:local}; those
through $\mu_0$ have sum at most $k_0+k_1+k_2=2$, so $k_1$ imposes no
condition. Once an even height is complete, the automaton only needs to
remember the pair $(\tau,\beta)$ read there; \S\ref{sec:primc} reduces
it to four states in this way. Figure~\ref{fig:primc-automata} shows
both automata for the boundary $(2,0,0)$.

\subsection{The completion and the verification theorem}\label{ss:completion}
Put
\[
 \A=\Q[q]\langle U,R\rangle/(RU-qUR),\qquad
 \widehat\A=\varprojlim_N \A/q^N\A=\Q[U,R][[q]],
 \qquad \widehat{\mathbb B}=\Q[R][[q]].
\]
The notation on the right refers to normal symbols; multiplication in
$\widehat\A$ remains the quantum product. Completion is taken
before inverting $q$, so the coefficient of each power of $q$ is a
polynomial in $U,R$, not an arbitrary series. For
$s=\sum_{k\ge0}q^ks_k(U,R)\in\widehat\A$, with $s_k$ in normal
form, the marginal and the evaluation of~\S\ref{sec:teaser} are
\[
 \rho(s)=\sum_{k\ge0}q^ks_k(1,R),\qquad
 \varepsilon(s)=\sum_{k\ge0}q^ks_k(1,1).
\]
Both are continuous and linear. A series $b=\sum_{n\ge0}b_nR^n$ with
coefficients in $\Q((q))$ belongs to $\widehat{\mathbb B}[q^{-1}]$
precisely when $\nu(b_n)\to+\infty$, and to $\widehat{\mathbb B}$ if
moreover all $b_n$ have nonnegative valuation. Similarly, a normally
ordered double series $\sum_{i,j}a_{i,j}U^iR^j$ belongs to
$\widehat\A[q^{-1}]$ precisely when $\nu(a_{i,j})\to+\infty$ as
$i+j\to\infty$. In practice we prove bounds of the form
\[
 \nu(a_{i,j})\ge\delta(i+j)-c,\qquad
 \nu(b_n)\ge\delta n-c\qquad(\delta>0).
\]

\begin{theorem}[Verification]\label{thm:verification}
Suppose $S\in\widehat\A$ and
\begin{equation}\label{eq:certificate-identity}
 S=(1-U)a+b(1-R)
\end{equation}
with $a,b\in\widehat\A[q^{-1}]$. Then $\varepsilon(S)=0$.
Conversely, if $\varepsilon(S)=0$, there is a unique pair
satisfying~\eqref{eq:certificate-identity} with
$a\in\widehat\A$ and $b\in\widehat{\mathbb B}$. If $D(x,R)$ is
the normal symbol of $S$, the normal symbols of this pair are
\begin{equation}\label{eq:canonical-boundaries}
 a(x,R)=\frac{D(x,R)-D(1,R)}{1-x},\qquad
 b(R)=\frac{D(1,R)}{1-R}.
\end{equation}
\end{theorem}
\begin{proof}
Left multiplication by $U$ and right multiplication by $R$ shift the
corresponding normal exponent without changing the coefficient.
Consequently $\varepsilon((1-U)a)=\varepsilon(b(1-R))=0$, first for
polynomials and then by continuity. This proves the first assertion.

Conversely, put $C=\rho(S)=D(1,R)$. At each $q$-degree, the polynomial
$D_k(x,R)-D_k(1,R)$ is divisible by $1-x$, and $C_k(1)=0$ implies that
$C_k(R)$ is divisible by $1-R$. These divisions
construct~\eqref{eq:canonical-boundaries} in the indicated completions,
and the symbol of $(1-U)a+b(1-R)$ is $(1-x)a+(1-R)b=D$. If $b$ is
independent of $U$, applying $\rho$ to~\eqref{eq:certificate-identity}
gives $C=b(1-R)$, which determines $b$; then $a$ is determined as well,
since left multiplication by $1-U$ is injective, coefficientwise in $q$.
\end{proof}

\subsection{Certificates}
Holonomicity is not needed for the telescoping argument itself; it
provides finite descriptions of $a$ and $b$. Set
\[
 \mathcal O=\K\langle R,\Ups\rangle/(\Ups R-qR\Ups),
\]
acting on $\K[[R]]$ by $\Ups g(R)=g(qR)$. An operator followed by a
polynomial in parentheses denotes the action on that polynomial; thus
$P(R)$ is the image of the series $R$ under $P\in\mathcal O$. A series
in $\K[[R]]$ is \emph{holonomic} if it is annihilated by a nonzero
element of $\mathcal O$, equivalently if its $q$-shift span over
$\K(R)$ is finite-dimensional. For a normal symbol in $x,R$ we require
the analogous finiteness of the joint shift span over $\K(x,R)$, under
$\Ups$ and $\Xi$, where $(\Xi g)(x,R)=g(qx,R)$.

A \emph{certificate} for a partition identity $\varepsilon(S)=0$
consists of
\begin{enumerate}[label=(c\arabic*)]
\item operator kernels defining $b$ over $\mathcal O$ and $a$ over its
      analogue in $x,R$;
\item initial conditions selecting the intended solutions of those
      kernels;
\item exact algebraic identities verifying~\eqref{eq:certificate-identity};
\item a proof that $\nu(b_n)\ge n+1$: a direct certificate
      (Definition~\ref{def:conv-cert}) or an independent proof of the
      identity (Corollary~\ref{cor:indirect}).
\end{enumerate}
No certificate is needed for the convergence of $a$, by
Theorem~\ref{thm:A}(2). The selecting data and~(c4) are essential: a
nonzero annihilating operator does not distinguish between its
solutions, and the coefficients of a formal solution need not tend
to~$0$.

\subsection{Constructing the holonomic kernels}\label{ss:kernels}
Consider a partition automaton. Let $\mathsf T(U,R)$ be the transition
matrix of its interior automaton, $r_0$ the row of weights of the
interior states reached by the first letter $R$, and $e_0$ the column of
weights of the interior states from which the last letter $U$ leads to
acceptance; thus $H=Rr_0(\I-\mathsf T(U,R))^{-1}e_0U$. The monomial
identity
\[
 RW(U,R)U=qUW(qU,qR)R,
\]
valid for every word $W$, gives a normalized realization
\begin{equation}\label{eq:normalized-realization}
 w\mathsf Q=r,\qquad H=we,\qquad
 \mathsf Q=\I-\mathsf T(qU,qR),\quad r=qUr_0,\quad e=e_0R.
\end{equation}
Every interior letter now carries a factor $q$, so
$\mathsf Q\equiv \I\pmod q$, and the selected solution is the geometric
series $w=r\sum_{k\ge0}(\I-\mathsf Q)^k$, which converges in
${\widehat\A}^{1\times s}$; this is why the word expansions
of~\S\ref{sec:teaser} converge. For a difference $S=H_1-H_2$ we use the
block-diagonal realization, with output column $(e_1,-e_2)$. More
generally, we allow entries of $\mathsf Q,r,e$ in
$\K\langle U,R\rangle/(RU-qUR)$ that are regular at $q=0$, provided
$\mathsf Q\equiv\I\pmod q$.

Define the anti-homomorphism $\Phi\colon\A\to\mathcal O$ by
\[
 \Phi(U^iR^j)=R^j\Ups^i.
\]
It expresses the right action on marginals: $\rho(uv)=\Phi(v)\rho(u)$
for normally ordered monomials. Hence
\begin{equation}\label{eq:marginal-system}
 \mathsf L_{ji}=\Phi(\mathsf Q_{ij}),\quad f=\rho(r)^{\top},\quad
 d_i=\Phi(e_i),\qquad \mathsf Lm=f,\quad C=\rho(S)=dm,
\end{equation}
where $m=\rho(w)^{\top}$. Since $\mathsf L\equiv \I\pmod q$, this
system has a unique solution in $\widehat{\mathbb B}^s$. Every
coefficient of $m$ and $C$ lies in $\K$: writing
$\mathsf L=\sum_jR^j\mathsf L_j(\Ups)$, the coefficient of $R^n$ is
obtained by solving with the matrix $\mathsf L_0(q^n)\equiv \I\pmod q$.
Since $C$ is linear in $f$, the factor $q$ of $f=q\,r_0^{\top}$ may be
moved into $d$; the examples of~\S\ref{sec:examples} do so.

Adjoin an auxiliary coordinate $u=1$, which satisfies $(\Ups-1)u=0$.
The canonical formal boundary $b$ is selected by the finite matrix
kernel
\begin{equation}\label{eq:b-kernel}
 \mathsf G_b\begin{pmatrix}m\\u\\b\end{pmatrix}=0,
 \qquad
 \mathsf G_b=\begin{pmatrix}
 \mathsf L&-f&0\\0&\Ups-1&0\\-d&0&1-R
 \end{pmatrix}.
\end{equation}
The last row means $(1-R)b=C$ in $\K[[R]]$, with no convergence
assumed. Thus, writing $C=\sum c_nR^n$,
\begin{equation}\label{eq:beta}
 b=\sum_{n\ge0}b_nR^n,\qquad b_n=\sum_{j=0}^nc_j.
\end{equation}
For path series as in Theorem~\ref{thm:A}, the boundary word of a
partition contains as many letters $R$ as its largest part. Hence $c_j$
is the difference of the weighted counts of the partitions with largest
part~$j$, and $b_n$ that of the partitions with largest part at
most~$n$.

For $a$, use the full right action on symbols,
$\Phi_x(U^iR^j)=x^iR^j\Ups^i$, and compile $\mathsf L_xm_x=f_x$, $D=d_xm_x$; a
subscript $1$ denotes $x=1$. The normal symbol $a$ is selected by
\begin{equation}\label{eq:a-kernel}
 \mathsf G_a\begin{pmatrix}m_x\\m_1\\u\\a\end{pmatrix}=0,\qquad
 \mathsf G_a=\begin{pmatrix}
 \mathsf L_x&0&-f_x&0\\
 0&\mathsf L_1&-f_1&0\\
 0&0&\Ups-1&0\\
 -d_x&d_1&0&1-x
 \end{pmatrix}.
\end{equation}
The last row is the first formula in~\eqref{eq:canonical-boundaries}.
These homogeneous kernels have sizes $s+2$ and $2s+2$; they are used
below to prove holonomicity, and require no scalar elimination. They do
not determine their solutions by themselves, since $(\Ups-1)u=0$ only
forces $u$ to be independent of $R$. Substituting $u=1$, that is,
deleting the auxiliary row and column, gives square inhomogeneous
systems of sizes $s+1$ and $2s+1$, with right-hand sides formed by $f$,
respectively $f_x$ and $f_1$. Their solutions in $\K[[R]]^{s+1}$ and
$\K[[x,R]]^{2s+1}$ are unique: the equations for $m$ and $m_1$ are
solved coefficientwise as above; those for $m_x$ are solved degree by
degree in $x,R$, since every letter carries a factor $x$ or $R$, so
that the part of $\I-\mathsf L_x$ that does not raise the total degree
comes from letter-free entries of $\I-\mathsf Q$ and is divisible by
$q$; and $1-R$, $1-x$ are injective.

\begin{lemma}[Finite elimination]\label{lem:finite-elimination}
Let $\mathsf G$ be a nonsingular $n\times n$ matrix over $\K(R)[\Ups]$,
of $\Ups$-degree at most $\delta$. For each column $k$, a finite linear
kernel computation gives a row $V$ and a nonzero scalar operator
$\Omega$ such that $V\mathsf G=\Omega e_k^{\top}$, with
$\deg_\Ups V\le(n-1)\delta$.
\end{lemma}
\begin{proof}
Put $t=(n-1)\delta$ and seek $V_j=\sum_{i=0}^tv_{ji}(R)\Ups^i$.
Annihilating the other columns gives at most $(n-1)(t+\delta+1)$ linear
equations over $\K(R)$ in $n(t+1)$ unknowns, one more than the number of
equations. A nonzero solution exists, and its remaining entry $\Omega$
cannot vanish because $\mathsf G$ is nonsingular. A common left scalar
denominator gives polynomial coefficients in $R$.
\end{proof}

The matrices~\eqref{eq:b-kernel} and~\eqref{eq:a-kernel} are
nonsingular over the corresponding Ore skew fields: normalization makes
each state matrix invertible in its $q$-adic operator completion, hence
nonsingular over the skew field (a left kernel vector, cleared of
denominators, would vanish in the completion), and the remaining
diagonal blocks $\Ups-1$, $1-R$, $1-x$ are nonzero. The lemma therefore proves that all
selected components are holonomic. The skew field is used only for
elimination; its inverses are never applied to series. For example,
$\Ups-1$ is invertible there, but annihilates the nonzero series $1$.

The symbol $a$ is jointly holonomic, not merely holonomic with $x$ as a
parameter. Left multiplication acts on symbols through the homomorphism
$U^iR^j\mapsto x^iR^j\Xi^j$; the opposite realization
$\mathsf Qv=e$, $S=rv$ therefore gives a pure $\Xi$-equation for $D$,
and the preceding construction a pure $\Ups$-equation. If their orders
are $p_1,p_2$, every mixed shift $\Xi^i\Ups^jD$ reduces to the span of
those with $i<p_1$, $j<p_2$: apply $\Xi^{i-p_1}\Ups^j$ to the first
equation to lower $i$ without increasing $j$, and $\Xi^i\Ups^{j-p_2}$ to
the second to lower $j$ without increasing $i$; shifted leading
coefficients remain nonzero in $\K(x,R)$. Finite shift rank is preserved
by sums and rational multiplication, and $C$ is holonomic in $R$ and
constant in $x$, so~\eqref{eq:canonical-boundaries} proves the assertion
for $a$.

The kernels can also be replaced by scalar operators with finite initial
data.

\begin{lemma}[Exceptional indices]\label{lem:exceptional}
Let $\Omega=\sum_{j\ge j_0}R^j\omega_j(\Ups)\in\mathcal O$ be nonzero,
with $\omega_{j_0}\ne0$, and put
$E_\Omega=\{n\ge0:\omega_{j_0}(q^n)=0\}$.
\begin{enumerate}
\item If $\Omega g=0$ for $g=\sum g_nR^n\in\K[[R]]$, then
\begin{equation}\label{eq:initial-recurrence}
 \omega_{j_0}(q^n)g_n
 =-\sum_{j>j_0}\omega_j(q^{n+j_0-j})g_{n+j_0-j},
\end{equation}
where negative subscripts mean zero. Hence $g$ is determined by its
coefficients at $E_\Omega$, and $g=0$ if and only if these vanish.
\item Let $\omega(Z)=\sum_i\gamma_iZ^i\in\K[Z]$ be nonzero, let $i_0$ be the
least index with $\gamma_{i_0}\ne0$, and let $n_0$ be the maximum of
$(\nu(\gamma_{i_0})-\nu(\gamma_i))/(i-i_0)$ over the indices $i>i_0$
with $\gamma_i\ne0$ ($n_0=-\infty$ if there are none). Then
$\nu(\omega(q^n))=\nu(\gamma_{i_0})+i_0n$ for every integer $n>n_0$; in
particular $\omega(q^n)\ne0$. Applied to $\omega=\omega_{j_0}$, this
shows that $E_\Omega$ is contained in $\{n\le n_0\}$, hence finite, and
$E_\Omega$ is computed by evaluating $\omega_{j_0}(q^n)\in\K$ for these
finitely many $n$.
\end{enumerate}
\end{lemma}
\begin{proof}
(1) is the coefficient of $R^{n+j_0}$ in $\Omega g=0$. (2) At $Z=q^n$,
the monomial $\gamma_iZ^i$ becomes $\gamma_iq^{in}$, of valuation
$\nu(\gamma_i)+in$. For $i>i_0$ with $\gamma_i\ne0$,
\[
 \nu(\gamma_{i_0})+i_0n<\nu(\gamma_i)+in
 \iff n>\frac{\nu(\gamma_{i_0})-\nu(\gamma_i)}{i-i_0}.
\]
Hence, for $n>n_0$, the term $\gamma_{i_0}q^{i_0n}$ has strictly
smaller valuation than every other term of $\omega(q^n)$, so the sum has
the valuation of this term.
\end{proof}

Values at $E_\Omega$ thus select at most one solution, and the
normalized construction supplies a consistent one.

\subsection{Convergence and proof of Theorem~\ref{thm:A}}
\label{ss:convergence}\label{ss:proof-A}
The kernels produce formal holonomic boundaries whether or not the
identity holds: the one-word partition language with series $S=RU=qUR$
gives $b=qR/(1-R)$, so $b_n=q$ for every $n\ge1$. For partition
languages, a single area estimate gives the precise criterion.

\begin{lemma}[Area bound]\label{lem:area}
In the setting of Theorem~\ref{thm:A}, write
$D(x,R)=\sum_{i,j}d_{i,j}x^iR^j$ and $C=D(1,R)=\sum_jc_jR^j$. Then
\begin{equation}\label{eq:partition-tail}
 \nu(d_{i,j})\ge i+j-1,\qquad
 \nu(c_j)\ge j,\qquad
 \nu\bigl(\varepsilon(S)-b_n\bigr)\ge n+1,\qquad
 \nu(a_{i,j})\ge i+j.
\end{equation}
\end{lemma}
\begin{proof}
A word in $R\{R,U\}^*U$ with $i$ letters $U$ and $j$ letters $R$ has
area at least $i+j-1$: count the pairs involving its first $R$ or its
last $U$. The weights have nonnegative valuation, so
$\nu(d_{i,j})\ge i+j-1$, and $\nu(c_j)\ge j$ because $i\ge1$. Next,
$\varepsilon(S)-b_n=\sum_{j>n}c_j$. Finally,
$(x^t-1)/(1-x)=-(1+x+\dots+x^{t-1})$ gives
$a_{i,j}=-\sum_{t>i}d_{t,j}$, of valuation at least
$\min_{t>i}(t+j-1)=i+j$.
\end{proof}

In particular, $(b_n)$ always converges to $\varepsilon(S)$, and
\begin{equation}\label{eq:convergence-equivalence}
 \begin{aligned}
 \varepsilon(S)=0
 &\quad\Longleftrightarrow\quad b_n\longrightarrow0\\
 &\quad\Longleftrightarrow\quad \nu(b_n)\ge n+1\text{ for all }n\\
 &\quad\Longleftrightarrow\quad b=q\Ups\phi,\quad\phi\in\Q[[q,R]].
 \end{aligned}
\end{equation}
Indeed, the last condition says exactly that
$\nu(q^{-n-1}b_n)\ge0$. Even arbitrarily slow decay along an unbounded
subsequence suffices: it forces $\varepsilon(S)=0$ and hence the
uniform bound $n+1$. Conversely, if $\varepsilon(S)\ne0$, then
$\nu(b_n)=\nu(\varepsilon(S))<n+1$ for all $n\ge\nu(\varepsilon(S))$,
a finite disproof.

\begin{corollary}[Indirect certificates]\label{cor:indirect}
If $\varepsilon(H_1)=\varepsilon(H_2)$ is proved by any means, then
the canonical pair $(a,b)$ of Theorem~\ref{thm:A}, selected by the
kernels of~\S\ref{ss:kernels}, satisfies $S=(1-U)a+b(1-R)$ with
$\nu(a_{i,j})\ge i+j$ and $\nu(b_n)\ge n+1$.
\end{corollary}
\begin{proof}
Combine Lemma~\ref{lem:area}, \eqref{eq:convergence-equivalence} and
Theorem~\ref{thm:verification}.
\end{proof}

\begin{proof}[Proof of Theorem~\ref{thm:A}]
For~(1), set $u=1$ in~\eqref{eq:a-kernel} and~\eqref{eq:b-kernel};
\S\ref{ss:kernels} proves uniqueness, that the coefficients $b_n$ lie
in $\K$, and holonomicity in the stated variables.
Lemma~\ref{lem:area} proves~(2), and~\eqref{eq:convergence-equivalence}
the equivalence of~(3a)--(3d). A kernel for $\phi=q^{-1}b(R/q)$ follows
by replacing $R$ by $R/q$ in a kernel for $b$. Condition~(3c)
implies~(3e), and Theorem~\ref{thm:verification}
gives~(3e)$\Rightarrow$(3a) together with the uniqueness statement.
\end{proof}

\subsection{Certifying convergence}\label{ss:conv-cert}
Items~(c1)--(c3) of a certificate are computed from the automata
(\S\ref{ss:kernels}), and by~\eqref{eq:convergence-equivalence}
item~(c4) is equivalent to the identity itself. We supply~(c4) in one of
two ways. A \emph{direct} certificate (Definition~\ref{def:conv-cert})
proves that the preimage $\phi=q^{-1}b(R/q)$ is $q$-integral by finitely
many exact checks on holonomic series. An \emph{indirect} certificate is
an independent proof of the identity, which
Corollary~\ref{cor:indirect} converts into convergence.

Let $\mathcal F\subset\K(U)$ be the ring of quotients of a polynomial in
$\Q[q,U]$ by a polynomial in $1+q\Q[q,U]$. For $g\in\mathcal F$, the
operator $g(\Ups)$ multiplies the coefficient of $R^n$ by
$g(q^n)\in\Q[[q]]$. Hence the \emph{integral operators}
$\sum_jR^jg_j(\Ups)$, with $g_j\in\mathcal F$, preserve $\Q[[q,R]]$;
they form a ring, since $g(\Ups)R^j=R^jg(q^j\Ups)$. Call a matrix
$\mathsf L$ of integral operators an \emph{integral kernel} if the
$R$-free part of $\I-\mathsf L$ has entries in $q\mathcal F$. Then
$\I-\mathsf L$ raises the total degree in $q$ and $R$, so for
$f\in\Q[[q,R]]^s$ the series $m=\sum_{k\ge0}(\I-\mathsf L)^kf$ converges
in $\Q[[q,R]]^s$ and solves $\mathsf Lm=f$. It is the only solution in
$\Q((q))[[R]]^s$: the coefficient of $R^n$ in $\mathsf Lm=f$ determines that
of $m$ from lower ones, through a matrix congruent to $\I$ modulo~$q$. If, in addition, $f\in\K[[R]]^s$, the same induction gives $m\in\K[[R]]^s$.

Integral kernels do not always suffice: in~\eqref{eq:euler-direct},
the boundary of Euler's identity combines two series whose coefficients
are $q$-adic units, and its convergence rests on a cancellation between
them. The following series account for such cancellations. For
$\ell\ge1$ and $\alpha\in\mathcal F$, the integral kernel
$1-R^\ell\alpha(\Ups)$ selects
\[
 u_\alpha=\sum_{k\ge0}\Bigl(\prod_{j=0}^{k-1}\alpha(q^{\ell j})\Bigr)
 R^{\ell k},\qquad\bigl(1-R^\ell\alpha(\Ups)\bigr)u_\alpha=1.
\]

\begin{proposition}[First-order series]\label{prop:dilation-checker}
Fix $\ell\ge1$. Let $\alpha,\beta$ be quotients of polynomials in
$1+q\Q[q,U]$, and suppose that $\theta=\beta/\alpha$ satisfies
$\theta-1\in q^cU^r\mathcal F$ for some integers $c\ge0$ and $r\ge1$.
\begin{enumerate}
\item The following are equivalent:
  \begin{enumerate}
  \item $\prod_{j\ge0}\theta(q^{\ell j})=1$;
  \item $u_\alpha-u_\beta=q^c\Ups^r\psi$ for some $\psi\in\Q[[q,R]]$;
  \item the coefficients of $u_\alpha-u_\beta$ tend to~$0$.
  \end{enumerate}
\item Condition~(a) holds if there are an integer $\lambda\ge2$ and a
  quotient $\chi$ of polynomials in $1+q\Q[q,U]$, with $\chi(0)=1$, such
  that
  \begin{equation}\label{eq:integral-dilation-check}
   \theta(U)\chi(U)=\prod_{i=0}^{\lambda-1}\chi(q^{\ell i}U^\lambda).
  \end{equation}
\end{enumerate}
\end{proposition}
\begin{proof}
The coefficient of $R^{\ell k}$ in $u_\alpha$ is
$u_k=\prod_{j<k}\alpha(q^{\ell j})\in1+q\Q[[q]]$, that in $u_\beta$ is
$u_k\pi_k$ with $\pi_k=\prod_{j<k}\theta(q^{\ell j})$, and all other
coefficients vanish. Since $\nu(\theta(q^{\ell j})-1)\ge c+r\ell j$, the
products $\pi_k$ converge to $\pi_\infty=\prod_{j\ge0}\theta(q^{\ell j})$.

(1) If $\pi_\infty=1$, then
$\pi_k=\prod_{j\ge k}\theta(q^{\ell j})^{-1}\equiv1\pmod{q^{c+r\ell k}}$,
so the coefficient $u_k(1-\pi_k)$ of $R^{\ell k}$ in $u_\alpha-u_\beta$
has valuation at least $c+r\ell k$. This is~(b), because
$[R^n]q^c\Ups^r\psi=q^{c+rn}[R^n]\psi$. Next, (b) implies~(c) because
$r\ge1$. Finally, (c) means that $u_k(1-\pi_k)\to0$; since the $u_k$ are
units, $\pi_k\to1$, so $\pi_\infty=1$.

(2) Since $\chi(0)=1$, we have $\chi-1\in U\mathcal F$, so each
$\chi(q^{\ell t})$ is a unit congruent to $1$ modulo $q^{\ell t}$.
Evaluate~\eqref{eq:integral-dilation-check} at $U=q^{\ell j}$ and
multiply over $0\le j<k$. As $(i,j)$ ranges over
$\{0,\dots,\lambda-1\}\times\{0,\dots,k-1\}$, the integer $i+\lambda j$
ranges over $\{0,\dots,\lambda k-1\}$, so
\[
 \pi_k=\frac{\prod_{t=0}^{\lambda k-1}\chi(q^{\ell t})}
            {\prod_{t=0}^{k-1}\chi(q^{\ell t})}
      =\prod_{t=k}^{\lambda k-1}\chi(q^{\ell t})\equiv1\pmod{q^{\ell k}}.
\]
Hence $\pi_\infty=1$.
\end{proof}

The coefficients of $u_\alpha$ and $u_\beta$ are $q$-adic units, so
neither series converges at $R=1$. Part~(1) reduces the convergence of
their difference to the product identity~(a), in which
$\theta(q^{\ell j})$ is the quotient of the ratios $\beta(q^{\ell j})$
and $\alpha(q^{\ell j})$ of consecutive coefficients of $u_\beta$ and
$u_\alpha$; in Corollary~\ref{cor:product} below, (a) is the partition
identity itself. Part~(2) certifies~(a) by a single identity
in~$\K(U)$, which we call a \emph{dilation witness}. First-order series arise as bounded counts of languages
defined by independent conditions on the multiplicities $f_h$
of~\eqref{eq:mult-word}.

\begin{corollary}[Periodic product languages]\label{cor:product}
Fix $\ell\ge1$ and, for each $i\in\Z/\ell$, eventually periodic sets
$M_{i,1},M_{i,2}\subseteq\Z_{\ge0}$ containing~$0$, with generating
functions $F_i(x)=\sum_{m\in M_{i,1}}x^m$ and
$G_i(x)=\sum_{m\in M_{i,2}}x^m$. Let $H_1$ be the path series of the
language of partitions in which, for every $h\ge1$, the multiplicity
$f_h$ of $h$ lies in $M_{h\bmod\ell,1}$, and define $H_2$ likewise with
the sets~$M_{i,2}$. Put
\[
 \alpha(U)=\prod_{i=1}^\ell F_i(q^iU),\qquad
 \beta(U)=\prod_{i=1}^\ell G_i(q^iU).
\]
If $\theta=\beta/\alpha$ has a dilation witness, then
$\varepsilon(H_1)=\varepsilon(H_2)$, and the canonical boundaries of
$H_1-H_2$ satisfy $\nu(a_{i,j})\ge i+j$ and $\nu(b_n)\ge n+1$.
\end{corollary}
\begin{proof}
The languages are regular by~\eqref{eq:mult-word}: an automaton needs
to remember only the height modulo~$\ell$ and the number of letters $U$
read at the current height, up to the sum of a preperiod and a period
of the sets $M_{i,1}$, $M_{i,2}$. Each $F_i(x)$ is a quotient of polynomials in
$1+x\Q[x]$, so $\alpha$, $\beta$ and $\theta$ are quotients of
polynomials in $1+q\Q[q,U]$; moreover $\theta(0)=1$, so
$\theta-1\in U\mathcal F$. By Proposition~\ref{prop:dilation-checker},
with $c=0$ and $r=1$, the coefficients of $u_\alpha-u_\beta$ tend
to~$0$. Grouping the part sizes $\ell j+1,\dots,\ell j+\ell$ for
$0\le j<k$, the partitions of the first language with largest part at
most $\ell k$ have generating function
$\prod_{h=1}^{\ell k}F_h(q^h)=\prod_{j<k}\alpha(q^{\ell j})$, and
similarly for the second. By the remark after~\eqref{eq:beta},
$b_{\ell k}$ is therefore the coefficient of $R^{\ell k}$ in
$u_\alpha-u_\beta$, so $b_{\ell k}\to0$. Since decay along a
subsequence suffices in~\eqref{eq:convergence-equivalence},
$\varepsilon(H_1)=\varepsilon(H_2)$ and $\nu(b_n)\ge n+1$; finally,
$\nu(a_{i,j})\ge i+j$ by Theorem~\ref{thm:A}(2).
\end{proof}

Glaisher's identity (Proposition~\ref{prop:glaisher}) is the basic
example.

\begin{definition}\label{def:conv-cert}
A \emph{direct certificate} for the convergence of a canonical
boundary~$b$ consists of
\begin{enumerate}
\item finitely many \emph{integral inputs}: components of solutions of
  $\mathsf Lm=f$, with $\mathsf L$ an integral kernel and
  $f\in\Q[q][R]^s$, and series $\psi$ as in
  Proposition~\ref{prop:dilation-checker}(1b), each with a dilation
  witness;
\item a series $\phi'$ obtained from the integral inputs and from
  $\Q[q][R]$ by sums and integral operators;
\item operator identities and selecting values proving $b=q\Ups\phi'$.
\end{enumerate}
\end{definition}

By construction $\phi'\in\Q[[q,R]]$, so a direct certificate proves
Theorem~\ref{thm:A}(3d), and hence all of Theorem~\ref{thm:A}(3). An
identity between selected solutions is checked either by explicit
operator multipliers or by a common annihilator and the finitely many
exceptional coefficients of Lemma~\ref{lem:exceptional}.

Euler's identity (Corollary~\ref{cor:euler-integral}) and the
five-state identity of~\S\ref{ss:nonlocal-memory} have direct
certificates. The other identities of~\S\ref{sec:examples} have indirect
ones: Glaisher's identity by Corollary~\ref{cor:product}, and those
of~\S\S\ref{sec:primc}--\ref{ss:berkovich-uncu} by finite telescoping
identities, a quadratic invariant (for Primc) and a recurrence with an
explicit solution (for Berkovich--Uncu), together with explicit
valuation estimates; the proof for Primc also uses the
G\"ollnitz--Gordon identities and Jacobi's triple product.
Appendix~\ref{app:searches} describes further finite tests, not needed
for these examples, which bound the valuations of $b_n$ instead of
exhibiting an integral preimage.

\subsection{Finite searches}\label{ss:searches}
The searches were developed to certify the vanishing differences
between the $222\,354$ minimal deterministic interior automata with at
most four states, enumerated with one initial state and partial
transitions: grouped by their series to order~$q^{50}$, they form
$19\,556$ non-singleton families. For bounded operator degrees, the
binding identities of a direct certificate, and the telescoping
certificates and quadratic invariants below, are solutions of finite
linear systems; first-order data and dilation witnesses of bounded
degree can be enumerated and checked
by~\eqref{eq:integral-dilation-check}; and rational shift quotients for
a proposed bounded sum are checked by clearing denominators, as
in~\S\ref{ss:berkovich-uncu}. Our \textsc{Julia} implementation (\S\ref{ss:implementation}) uses
such bounded searches, and Appendix~\ref{app:searches} describes further
ones.

\subsubsection*{Telescoping}
Let $(f_n)_{n\ge0}$ satisfy $f_{n+1}=r(q^n)f_n$ for a rational function
$r$, in which further parameters, such as $v$ in~\S\ref{sec:primc}, may
occur. A \emph{telescoping certificate} for a residual $\varpi$ is a
rational function $\Psi$ such that
\begin{equation}\label{eq:telescoping}
 r(\zeta)\Psi(q\zeta)-\Psi(\zeta)=\varpi(\zeta).
\end{equation}
This is an identity between rational functions of an indeterminate
$\zeta$, checked by clearing denominators; for a fixed denominator and a
degree bound, $\Psi$ is found by solving a linear system. Evaluating at
$\zeta=q^n$, multiplying by $f_n$ and summing over $0\le n\le N$ gives
\[
 \sum_{n=0}^Nf_n\varpi(q^n)=f_{N+1}\Psi(q^{N+1})-f_0\Psi(1),
\]
provided that no denominator vanishes at the points $q^n$; by
Lemma~\ref{lem:exceptional}(2), this takes finitely many evaluations.
For an infinite sum, the upper flux $f_{N+1}\Psi(q^{N+1})$ must also
tend to~$0$. This is checked symbolically, treating $N$ and $q^N$ as
independent symbols. By Lemma~\ref{lem:exceptional}(2), applied to
numerators and denominators, $\nu(r(q^n))$ and $\nu(\Psi(q^n))$ are
affine functions of $n$ for large $n$, read off from the terms of lowest
degree in $\zeta$. Hence
\[
 \nu\bigl(f_{N+1}\Psi(q^{N+1})\bigr)
 =\nu(f_0)+\sum_{n=0}^N\nu\bigl(r(q^n)\bigr)+\nu\bigl(\Psi(q^{N+1})\bigr)
\]
is a quadratic polynomial in $N$ for large $N$, and the flux tends to
$0$ when this polynomial tends to $+\infty$, typically because it is at
least $c'N^2-c$ with $c'>0$. The same count applies to other valuations,
such as the $v$-adic one in~\S\ref{sec:primc}. All checks are thus
manipulations of rational functions in $\zeta$ and of polynomials
in~$N$, which any computer algebra system performs. Several summation
variables are handled by sums of such differences.

\subsubsection*{Quadratic invariants}
For a system $F(z)=\mathsf M(z)F(q^hz)$ with $F$ formal at $z=0$, a
\emph{quadratic invariant} is a polynomial symmetric matrix
$\mathsf J(z)$ such that
\[
 \mathsf M(z)^{\top}\mathsf J(z)\mathsf M(z)=\lambda(z)\mathsf J(q^hz)
\]
for some polynomial $\lambda(z)$. For fixed $\lambda$, this equation is
linear in the coefficients of $\mathsf J$. If $\lambda(0)=0$, it forces
$F^{\top}\mathsf JF=0$, because this series $g(z)$ satisfies
$g(z)=\lambda(z)g(q^hz)$ and so has no term of lowest order. Such a
relation can reveal a factorization. In~\S\ref{sec:primc}, it is
$\det\mathsf H=0$ for a $2\times2$ matrix $\mathsf H$ formed from the
four states, on which the transfer acts by
$\mathsf H(v)=\mathsf\Gamma_+(v)\mathsf H(qv)\mathsf\Gamma_-(v)^{\top}$,
that is, through the tensor product of two two-state transfers. The
solution is the rank-one matrix $u(v)u(-v)^{\top}$, and the components
of $u$ are found by telescoping.

\section{Partition identities}\label{sec:examples}

\subsection{Glaisher's identity}\label{ss:glaisher}
For an integer $\lambda\ge2$, let $H_1$ and $H_2$ be the path series of
the languages of partitions with no part divisible by $\lambda$,
respectively in which no part occurs $\lambda$ or more times. They are
as in Corollary~\ref{cor:product}, with $\ell=\lambda$,
$M_{0,1}=\{0\}$, $M_{i,1}=\Z_{\ge0}$ for $i\not\equiv0$, and
$M_{i,2}=\{0,\dots,\lambda-1\}$; thus
\[
 F_i(x)=\frac1{1-x}\quad(i\not\equiv0),\qquad F_0(x)=1,\qquad
 G_i(x)=1+x+\dots+x^{\lambda-1}.
\]
For $\lambda=2$, they are the series $H_{\mathrm o}$ and
$H_{\mathrm d}$ of~\S\ref{sec:teaser}.

\begin{proposition}[following Glaisher~\cite{Glaisher1883}]\label{prop:glaisher}
Let $\lambda\ge2$ be an integer. For every $n$, the partitions of $n$
with no part divisible by $\lambda$ are as many as those in which no part occurs
$\lambda$ or more times; for $\lambda=2$, this is Euler's identity. The
canonical boundaries of $H_1-H_2$ satisfy $\nu(a_{i,j})\ge i+j$ and
$\nu(b_n)\ge n+1$.
\end{proposition}
\begin{proof}
In the notation of Corollary~\ref{cor:product}, since
$G_i(x)=(1-x^\lambda)/(1-x)$,
\[
 \alpha(U)=\prod_{i=1}^{\lambda-1}\frac1{1-q^iU},\qquad
 \beta(U)=\prod_{i=1}^{\lambda}\frac{1-q^{\lambda i}U^\lambda}{1-q^iU}.
\]
The factors $1-q^iU$ with $i<\lambda$ cancel in $\theta=\beta/\alpha$, so
\[
 \theta(U)=\frac{\prod_{i=1}^{\lambda}(1-q^{\lambda i}U^\lambda)}
                {1-q^\lambda U}.
\]
The polynomial $\chi(U)=1-q^\lambda U$ is a dilation witness: both sides
of~\eqref{eq:integral-dilation-check} equal the numerator of $\theta$.
Corollary~\ref{cor:product} concludes.
\end{proof}

The dilation witness encodes the classical argument
$1+x+\dots+x^{\lambda-1}=(1-x^\lambda)/(1-x)$. The proof identifies
$b_{\lambda k}$ through its combinatorial meaning. For $\lambda=2$, the
same first-order data give the direct certificate of the teaser
(Definition~\ref{def:conv-cert}), whose checks use only the kernels.

\begin{corollary}[Euler's certificate]\label{cor:euler-integral}
In the notation of~\S\ref{sec:teaser}, the canonical boundary $b$ of
$H_{\mathrm o}-H_{\mathrm d}$ satisfies~\eqref{eq:euler-direct}, and
$b=q\Ups b'$ with $b'=qw+R(1-q\Ups)^{-1}(w+\Ups v)\in\Q[[q,R]]$.
\end{corollary}
\begin{proof}
For $\lambda=2$, the data in the proof of
Proposition~\ref{prop:glaisher} are
\[
 \alpha(U)=\frac1{1-qU},\qquad \beta(U)=B(U)=(1+qU)(1+q^2U),\qquad
 \chi(U)=1-q^2U,
\]
and $\theta(U)=(1-qU)(1+qU)(1+q^2U)$. Thus $u=u_\alpha$ and $v=u_\beta$
are the series of~\S\ref{sec:teaser}, and
Proposition~\ref{prop:dilation-checker}, with $\ell=2$, $c=2$ and
$r=1$, gives again $u-v=q^2\Ups w$ with $w\in\Q[[q,R]]$. In the
computations below, $g(\Ups)R^\ell=R^\ell g(q^\ell\Ups)$.

Put $B_{\mathrm o}=(1+R(1-q\Ups)^{-1})u$ and
$B_{\mathrm d}=(1+R(1+q\Ups))v$. By the equations of $u$ and $v$,
\[
 \begin{aligned}
 (1-q\Ups-q^2R^2)\,\Ups(1-q\Ups)^{-1}u
 &=\Ups\bigl(u-R^2(1-q\Ups)^{-1}u\bigr)=1,\\
 (1-qR-q^2R\Ups)\,\Ups B_{\mathrm d}
 &=\Ups\bigl(1-R(1+q\Ups)\bigr)B_{\mathrm d}
 =\Ups\bigl(1-R^2B(\Ups)\bigr)v=1.
 \end{aligned}
\]
By uniqueness of normalized solutions,
$m_{\mathrm o,1}=\Ups(1-q\Ups)^{-1}u$ and
$m_{\mathrm d,1}=\Ups B_{\mathrm d}$. The same equations give
\[
 \begin{aligned}
 (1-R)B_{\mathrm o}&=1+R\bigl((1-q\Ups)^{-1}-1\bigr)u=1+qRm_{\mathrm o,1},\\
 (1-R)B_{\mathrm d}&=1+qR\Ups B_{\mathrm d}=1+qRm_{\mathrm d,1}.
 \end{aligned}
\]
Hence $(1-R)(B_{\mathrm o}-B_{\mathrm d})=D(1,R)$, so
$b=B_{\mathrm o}-B_{\mathrm d}$, which is~\eqref{eq:euler-direct}.
Finally, $u=v+q^2\Ups w$,
$(1-q\Ups)^{-1}-(1+q\Ups)=q^2\Ups^2(1-q\Ups)^{-1}$ and $q^2R\Ups=q\Ups R$
give
\[
 B_{\mathrm o}-B_{\mathrm d}
 =q^2\Ups w+q^2R\Ups(1-q\Ups)^{-1}(w+\Ups v)=q\Ups b'.\qedhere
\]
\end{proof}

\subsection{Primc's three-row identities at level two}\label{sec:primc}
We now prove Theorem~\ref{p2:thm:main}, following the strategy
of~\S\S\ref{sec:teaser}--\ref{sec:certificates}. Both partition classes
are recognized by automata (Figure~\ref{fig:primc-automata}). In~$\A$,
left multiplication turns the automaton of the admissible partitions
into a linear system of $q$-difference equations in the variable~$x$
of the normal symbols, which marks the number of parts; at $R=1$, it
reduces to a four-state system. A quadratic invariant factors this
system into two two-state systems, and two telescoping certificates
solve these in closed form; both were found by the searches
of~\S\ref{ss:searches}. At $x=1$, the G\"ollnitz--Gordon identities and
Jacobi's triple product then give the products of
Theorem~\ref{p2:thm:main}, and Corollary~\ref{cor:indirect} converts
the resulting equalities into convergent holonomic certificates for the
automata.

\subsubsection{From automata to four states}
The admissible partitions are recognized by the window automaton of
Lemma~\ref{lem:window}, with $r=W=3$, $M=2$ and $K=2$, and the residue
side by the congruence automaton with $M=8$ (\S\ref{ss:compile});
Figure~\ref{fig:primc-automata} shows both for the boundary $(2,0,0)$.
We use the multiplicities $\tau_n,\beta_n,\mu_n$ and the boundary values
defined before Theorem~\ref{p2:thm:main}. Once an even height is
complete, the window automaton only needs to remember the pair
$(\tau,\beta)$ read there, and the boundary only fixes the initial pair
$(\tau_0,\beta_0)=(k_2,k_0)$. For a pair $(\tau,\beta)$, let
$F_{\tau,\beta}(z)$ be the generating function of the \emph{admissible
tails}: the finitely supported sequences $(\mu_n,\tau_n,\beta_n)_{n\ge1}$
satisfying~\eqref{p2:eq:local} with $(\tau_0,\beta_0)=(\tau,\beta)$,
with $z$ marking the number of parts and $q$ their total size, measured
from the even height at which $(\tau,\beta)$ was read; the zero tail
contributes~$1$. The next odd and even heights carry $\mu$ parts of
size~$1$ and $\tau'+\beta'$ parts of size~$2$, and shift every later
part by~$2$; this gives the factor $(zq)^\mu(zq^2)^{\tau'+\beta'}$ and
the substitution $z\mapsto q^2z$, while~\eqref{p2:eq:local} allows
$0\le\mu\le2-\max(\tau,\tau')-\max(\beta,\beta')$. Hence
\[
 F_{\tau,\beta}(z)=\sum_{\tau',\beta'\ge0}(zq^2)^{\tau'+\beta'}
 \sum_{\mu=0}^{2-\max(\tau,\tau')-\max(\beta,\beta')}(zq)^\mu
 F_{\tau',\beta'}(q^2z),
\]
where a negative upper bound contributes zero. This is the opposite
realization of~\S\ref{ss:kernels}, before normalization, with the
letters of two consecutive heights grouped together. Since
$F_{\tau,\beta}=F_{\beta,\tau}$ by top--bottom reflection, the four
states $00$, $01$, $02$, $11$ remain. For instance, after the pair
$(0,0)$, the next pair $(0,0)$ allows $\mu\le2$, the pairs $(0,1)$ and
$(1,0)$ allow $\mu\le1$, and the pairs $(0,2)$, $(2,0)$ and $(1,1)$ allow
only $\mu=0$; this is the first row of $\mathsf M$ below. Set
$F=(F_{00},F_{01},F_{02},F_{11})^{\top}$, $s=qz$ and $t=q^2z$. Then
\begin{equation}\label{p2:eq:transfer}
 F(z)=\mathsf M(z)F(q^2z),\qquad
 \mathsf M(z)=\begin{pmatrix}
 1+s+s^2&2t(1+s)&2t^2&t^2\\
 1+s&t(2+s)&t^2&t^2\\
 1&t&t^2&0\\
 1&2t&0&t^2
 \end{pmatrix}.
\end{equation}
The empty tail gives $F(0)=(1,1,1,1)^{\top}$. Since $\mathsf M(0)$ is an
idempotent of rank one, $\I-q^{2n}\mathsf M(0)$ is invertible for $n\ge1$, and
coefficient induction proves uniqueness in $\Q(q)[[z]]^4$. The four
generating functions of Theorem~\ref{p2:thm:main} are, in order,
$\Pi_{2,0,0}=F_{02}(1)$, $\Pi_{1,1,0}=F_{01}(1)$, $\Pi_{1,0,1}=F_{11}(1)$
and $\Pi_{0,2,0}=F_{00}(1)$. With $z$ marking the number of parts,
$F_{02},F_{00},F_{01},F_{11}$ are the series
$D^{(2)}_{2\Lambda_0}$, $D^{(2)}_{2\Lambda_1}$, $D^{(2)}_{\Lambda_0+\Lambda_1}$,
$D^{(2)}_{\Lambda_0+\Lambda_2}$ of~\cite[\S4]{KRTW}.

\subsubsection{Factoring the four-state system}
The search for quadratic invariants of~\S\ref{ss:searches}, applied
to~\eqref{p2:eq:transfer}, finds
\[
 \mathsf J(z)=\begin{pmatrix}0&0&0&1\\0&-2&0&0\\0&0&2qz&0\\1&0&0&0\end{pmatrix},
 \qquad \lambda(z)=-q^5z^3(1-qz),
\]
and polynomial arithmetic checks that
$\mathsf M(z)^{\top}\mathsf J(z)\mathsf M(z)=\lambda(z)\mathsf J(q^2z)$.
  The search takes $\mathsf J$ symmetric of degree at most one in~$z$,
  with $20$ unknown coefficients; comparing determinants, with
  $\det\mathsf M(z)=q^{11}z^6(1-qz)^2$, suggests
  $\lambda=\pm q^5z^3(1-qz)$. The ten entries of
  $\mathsf M^{\top}\mathsf J\mathsf M-\lambda\mathsf J(q^2z)$ have degree
  at most five in~$z$, so this is a $60\times20$ linear system over
  $\Q(q)$; its kernel is spanned by $\mathsf J$ for the minus sign and is
  zero for the plus sign. Since $\lambda(0)=0$, this
forces $F^{\top}\mathsf JF=0$, that is,
\begin{equation}\label{p2:eq:quadratic}
 F_{00}(z)F_{11}(z)=F_{01}(z)^2-qzF_{02}(z)^2.
\end{equation}
This relation suggests assembling the four series into a matrix of
determinant zero. Introduce $v$ by $v^2=qz$ and define
\[
 \mathsf H(v)=\begin{pmatrix}
 F_{00}(v^2/q)&F_{01}(v^2/q)+vF_{02}(v^2/q)\\
 F_{01}(v^2/q)-vF_{02}(v^2/q)&F_{11}(v^2/q)
 \end{pmatrix}
\]
and
\[
 \mathsf\Gamma_\pm(v)=\begin{pmatrix}1\pm v+v^2&qv^2\\1&qv^2\end{pmatrix}.
\]
By~\eqref{p2:eq:transfer}, direct multiplication gives the exact
identity
\begin{equation}\label{p2:eq:matrix-factor}
 \mathsf H(v)=\mathsf\Gamma_+(v)\mathsf H(qv)\mathsf\Gamma_-(v)^{\top};
\end{equation}
thus the transfer acts on $\mathsf H$ through the tensor product of two
two-state transfers, and~\eqref{p2:eq:quadratic} reads
$\det\mathsf H=0$. Equation~\eqref{p2:eq:matrix-factor} has a unique
solution with the all-ones constant matrix.
Indeed, put $\mathsf\Gamma_0=\mathsf\Gamma_\pm(0)$ and
$\PP(\mathsf Z)=\mathsf\Gamma_0\mathsf Z\mathsf\Gamma_0^{\top}$ for
$2\times2$ matrices $\mathsf Z$, and write
$\mathsf H(v)=\sum_n\mathsf H_nv^n$. The coefficient of $v^n$
in~\eqref{p2:eq:matrix-factor} determines $(\I-q^n\PP)\mathsf H_n$ from
$\mathsf H_0,\dots,\mathsf H_{n-1}$; since $\PP^2=\PP$, the
operator $\I-q^n\PP$ has, for $n>0$, the explicit inverse
\[
 (\I-q^n\PP)^{-1}=\I+\frac{q^n}{1-q^n}\PP,
\]
and for $n=0$ the image of $\PP$ consists of the multiples of the
all-ones matrix.

Let $u(v)=(g(v),h(v))^{\top}$ solve
\begin{equation}\label{p2:eq:spinor}
 u(v)=\mathsf\Gamma_+(v)u(qv),\qquad u(0)=(1,1)^{\top}.
\end{equation}
The same argument proves uniqueness. Since $\mathsf\Gamma_-(v)=\mathsf\Gamma_+(-v)$,
the matrix $u(v)u(-v)^{\top}$ solves \eqref{p2:eq:matrix-factor}
with the correct constant matrix. Consequently
\begin{align}
 F_{00}(v^2/q)&=g(v)g(-v),&
 F_{11}(v^2/q)&=h(v)h(-v),\label{p2:eq:diagonal}\\
 F_{01}(v^2/q)&=\frac{g(v)h(-v)+h(v)g(-v)}2,&
 F_{02}(v^2/q)&=\frac{g(v)h(-v)-h(v)g(-v)}{2v}.
 \label{p2:eq:cross}
\end{align}

\subsubsection{Two finite telescoping certificates}
The solution of \eqref{p2:eq:spinor} is
\begin{equation}\label{p2:eq:single-sums}
 g(v)=\sum_{n\ge0}\frac{q^{n(n-1)/2}v^n(-v;q)_n}{(q;q)_n},\qquad
 h(v)=\sum_{n\ge0}\frac{q^{n(n+1)/2}v^n(-v;q)_n}{(q;q)_n}.
\end{equation}
To verify this, set
\[
 f_n(v)=\frac{q^{n(n-1)/2}v^n(-v;q)_n}{(q;q)_n},\qquad \zeta=q^n,
\]
so that $h=\sum_n\zeta f_n$ and
\[
 \frac{f_{n+1}(v)}{f_n(v)}=r(\zeta):=\frac{v\zeta(1+v\zeta)}{1-q\zeta},\qquad
 \frac{f_n(qv)}{f_n(v)}=\frac{\zeta(1+v\zeta)}{1+v}.
\]
The two telescoping certificates
\[
 \Psi_0(\zeta)=\zeta-1,\qquad \Psi_1(\zeta)=-\frac{v\zeta(1-\zeta)}{1+v}
\]
satisfy~\eqref{eq:telescoping} with the following residuals:
\begin{align}
 r(\zeta)\Psi_0(q\zeta)-\Psi_0(\zeta)&=1-\zeta-v\zeta(1+v\zeta),\label{p2:eq:wz0}\\
 r(\zeta)\Psi_1(q\zeta)-\Psi_1(\zeta)&=
 \frac{v\zeta(1-\zeta)-qv^2\zeta^2(1+v\zeta)}{1+v}.
 \label{p2:eq:wz1}
\end{align}
Multiply by $f_n(v)$ and sum from $0$ to $N$; the left sides telescope.
Both lower fluxes vanish because $\Psi_i(1)=0$. The upper fluxes are
$f_{N+1}(v)\Psi_i(q^{N+1})$, of $v$-valuation at least $N+1$, so they
tend to zero in $\Q(q)[[v]]$. Equations
\eqref{p2:eq:wz0}--\eqref{p2:eq:wz1} therefore give
\[
 g(v)-h(v)=v(1+v)g(qv),\qquad
 h(v)=g(qv)+qv^2h(qv).
\]
These are precisely~\eqref{p2:eq:spinor}, and uniqueness
proves~\eqref{p2:eq:single-sums}. Both certificates were found by the
search of~\S\ref{ss:searches}: with $\Psi$ a polynomial in $\zeta$ of
  degree at most one, respectively two, and an unknown scalar multiple of
  the residual, clearing the denominator $1-q\zeta$ of $r$
  turns~\eqref{eq:telescoping} into linear systems over $\Q(q,v)$ of sizes
  $4\times3$ and $5\times4$, each with one-dimensional kernel.
  
\subsubsection{G\"ollnitz--Gordon specialization}
Put $q=r^2$ and $v=r$. The two sums in \eqref{p2:eq:single-sums} become
the classical G\"ollnitz--Gordon series:
\begin{align}
 g(r;r^2)&=\sum_{n\ge0}\frac{r^{n^2}(-r;r^2)_n}{(r^2;r^2)_n}
 =\frac1{(r,r^4,r^7;r^8)_\infty},\label{p2:eq:gg1}\\
 h(r;r^2)&=\sum_{n\ge0}\frac{r^{n^2+2n}(-r;r^2)_n}{(r^2;r^2)_n}
 =\frac1{(r^3,r^4,r^5;r^8)_\infty}.
 \label{p2:eq:gg2}
\end{align}
The final equalities are the G\"ollnitz--Gordon
identities~\cite{Gollnitz1967,Gordon1965}. Since $n^2$ and $n^2+2n$ have the parity of $n$, replacing $r$
by $-r$ in~\eqref{p2:eq:gg1}--\eqref{p2:eq:gg2} gives the products for
$g(-r;r^2)$ and $h(-r;r^2)$.

The series $g$, $h$ and $F_{\tau,\beta}(v^2/q)$ lie in $\Q[[q,v]]$ (a
tail with $k$ parts has size at least $k$, and $z^k=q^{-k}v^{2k}$), and
$q\mapsto r^2$, $v\mapsto r$ defines a continuous homomorphism
$\Q[[q,v]]\to\Q[[r]]$; so the identities above survive this
specialization, which corresponds to $z=1$.

Pairing the odd factors in the products at $r$ and $-r$ gives
\begin{align*}
 F_{00}(1;q)&=g(r;r^2)g(-r;r^2)
 =\frac1{(q,q^7;q^8)_\infty(q^2;q^4)_\infty^2},\\
 F_{11}(1;q)&=h(r;r^2)h(-r;r^2)
 =\frac1{(q^3,q^5;q^8)_\infty(q^2;q^4)_\infty^2}.
\end{align*}
These are the products for the last and third rows of
Theorem~\ref{p2:thm:main}: residues $1,7$, respectively $3,5$, have one
colour, and residues $2,6$ have two colours.

\subsubsection{The two cross-products}
Write $\theta(y;p)=(y,p/y;p)_\infty$. Jacobi's triple
product~\cite[Theorem~2.8]{Andrews} gives, in $\Z[y_1^{\pm1},y_2^{\pm1}][[p]]$,
\begin{align}
 \theta(-y_1;p)\theta(y_2;p)+\theta(y_1;p)\theta(-y_2;p)
 &=2\frac{(p^2;p^2)_\infty^2}{(p;p)_\infty^2}
 \theta(y_1y_2;p^2)\theta(py_1/y_2;p^2),\label{p2:eq:theta-plus}\\
 \theta(-y_1;p)\theta(y_2;p)-\theta(y_1;p)\theta(-y_2;p)
 &=2y_1\frac{(p^2;p^2)_\infty^2}{(p;p)_\infty^2}
 \theta(y_2/y_1;p^2)\theta(py_1y_2;p^2).
 \label{p2:eq:theta-minus}
\end{align}
These two dissections follow by multiplying the left sides by
$(p;p)_\infty^2$ and expanding the two bilateral series; the sign factor
is $(-1)^n\pm(-1)^m$. For the plus sign put $m=k+l$, $n=k-l$: the
exponent becomes $k^2+l^2-k$ and the monomial $(y_1y_2)^k(y_1/y_2)^l$.
For the minus sign put $m=k+l+1$, $n=k-l$: the exponent becomes
$k^2+l^2+l$ and the monomial $y_1(y_1y_2)^k(y_1/y_2)^l$. In each case
the two independent sums are triple products.

The products $g(r)h(-r)$ and $h(r)g(-r)$ have common denominator
\[
 \mathcal D=(r^4;r^8)_\infty^2
 (r^2,r^6,r^{10},r^{14};r^{16})_\infty,
\]
and their sum and difference have numerators
$\theta(-r;r^8)\theta(r^3;r^8)\pm\theta(r;r^8)\theta(-r^3;r^8)$.
Take $y_1=r$, $y_2=r^3$, $p=r^8$ in the two dissections. Since
$(r^{16};r^{16})_\infty^2/(r^8;r^8)_\infty^2=(r^8;r^{16})_\infty^{-2}$,
these numerators are
\[
 \frac{2(r^4,r^6,r^{10},r^{12};r^{16})_\infty}
 {(r^8;r^{16})_\infty^2},\qquad
 \frac{2r(r^2,r^4,r^{12},r^{14};r^{16})_\infty}
 {(r^8;r^{16})_\infty^2}.
\]
Substituting into \eqref{p2:eq:cross} and cancelling factors gives
\begin{align*}
 F_{01}(1;r^2)&=
 \frac1{(r^2,r^4,r^8,r^8,r^{12},r^{14};r^{16})_\infty},\\
 F_{02}(1;r^2)&=
 \frac1{(r^4,r^6,r^8,r^8,r^{10},r^{12};r^{16})_\infty},
\end{align*}
giving the second and first rows of Theorem~\ref{p2:thm:main}.
In both products residue $4$ has two colours; the other four listed
residues have one colour each. The bilateral series converge
$r$-adically, since the triple-product summands have positive quadratic
valuation in both summation indices.
This completes the proof of Theorem~\ref{p2:thm:main}.

\subsubsection{Certificates for the automata}
Let $S$ be the difference of the two automata of
Figure~\ref{fig:primc-automata}. For the boundary $(2,0,0)$, they have $29$ and $9$ interior states, so $s=38$, and the homogeneous kernels
$\mathsf G_a,\mathsf G_b$ of~\eqref{eq:a-kernel}
and~\eqref{eq:b-kernel} have sizes $78$ and $40$. The proof above
establishes $\varepsilon(S)=0$, so Corollary~\ref{cor:indirect}
supplies convergent holonomic certificates for the original automata,
with $\nu(a_{i,j})\ge i+j$ and $\nu(b_n)\ge n+1$, without eliminating
their matrices to scalar operators. The implementation
(\S\ref{ss:implementation}) lists both automata and computes the sparse
operator entries of these kernels.

\subsection{The Rogers--Ramanujan identities}\label{ss:rr}
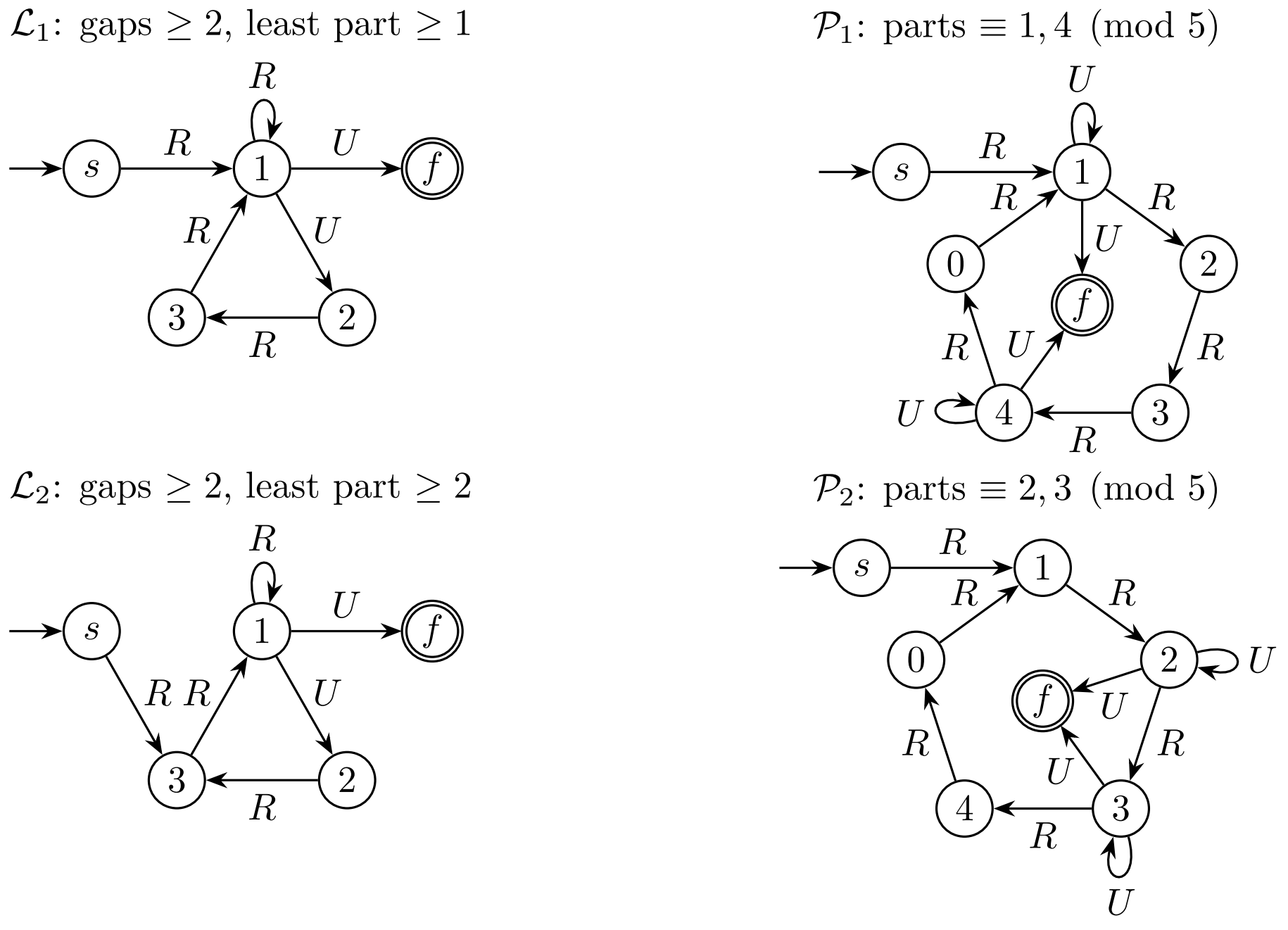
\begin{figure}[htbp]
\centering
\tikzset{
  rr/.style={
    >=Stealth, semithick,
    initial text={}, initial distance=5mm,
    every state/.style={
      draw,circle,minimum size=15pt,inner sep=1pt
    },
    every edge/.style={draw,->},
    every loop/.style={looseness=7}
  },
  lab/.style={
    postaction={decorate},
    decoration={
      markings,
      mark=at position .5 with {
        \node[
          anchor=center,inner sep=0pt,outer sep=0pt,
          transform shape=false
        ] at (0,2.5mm) {$#1$};
      }
    }
  }
}
\foreach \i/\a/\b in {1/1/4,2/2/3}{%
  \begin{minipage}[t]{.46\linewidth}
  \centering\vspace{0pt}
  $\mathcal L_{\i}$: gaps $\ge2$, least part $\ge\i$
  \par\medskip
  \begin{tikzpicture}[rr]
    \path[use as bounding box]
      (-2.4,-1.8) rectangle (2,1);
    \node[state,initial] (s) at (-1.6,0) {$s$};
    \node[state] (g1) at (0,0) {$1$};
    \node[state] (g2) at (0.8,-1.4) {$2$};
    \node[state] (g3) at (-0.8,-1.4) {$3$};
    \node[state,accepting] (f) at (1.6,0) {$f$};

    \ifnum\i=1
      \path (s) edge[lab=R] (g1);
    \else
      \path (s) edge[lab=R] (g3);
    \fi

    \path
      (g1) edge[loop above,lab=R] (g1)
           edge[lab=U] (g2)
           edge[lab=U] (f)
      (g2) edge[lab=R] (g3)
      (g3) edge[lab=R] (g1);
  \end{tikzpicture}
  \end{minipage}\hfill
  \begin{minipage}[t]{.50\linewidth}
  \centering\vspace{0pt}
  $\mathcal P_{\i}$: parts $\equiv\a,\b\pmod5$
  \par\medskip
  \begin{tikzpicture}[rr]
    \foreach \r/\angle in
      {0/162,1/90,2/18,3/-54,4/-126}
      \node[state] (r\r) at (\angle:1.25) {$\r$};

    \node[state,initial] (s) at (-1.7,1.25) {$s$};
    \node[state,accepting] (f) at (0,0) {$f$};

    \path (s) edge[lab=R] (r1);

    \foreach \j/\k in {0/1,1/2,2/3,3/4,4/0}
      \path (r\j) edge[lab=R] (r\k);

    \foreach \r in {\a,\b}
      \path (r\r) edge[lab=U] (f);

    \ifnum\i=1
      \path
        (r1) edge[loop above,lab=U] (r1)
        (r4) edge[loop left,lab=U] (r4);
    \else
      \path
        (r2) edge[loop right,lab=U] (r2)
        (r3) edge[loop below,lab=U] (r3);
    \fi
  \end{tikzpicture}
  \end{minipage}
  \par\medskip
}
\caption{Automata for the Rogers--Ramanujan identities:
gap conditions (left) and residue conditions (right).
Each automaton reads the complete boundary word, from its unique
initial state $s$ to its unique final state $f$.
Numbers on the right record height modulo $5$.}
\label{fig:rr-automata}
\end{figure}

This classical case shows that the whole proof, product side included, reduces to finite rational identities of the kind found by the searches of~\S\ref{ss:searches}. Write $(q)_n=(q;q)_n$. For $i=1,2$, the language
\[
 \mathcal L_i=R^iR^*U\,(R^2R^*U)^*
\]
encodes nonempty partitions with successive parts differing by at least
two and smallest part at least $i$: successive blocks of $R$'s record
the smallest part and the successive gaps. Subtracting the staircase
$i,i+2,\ldots,i+2(n-1)$ gives
\[
 1+\varepsilon(H_{\mathcal L_i})
 =\sum_{n\ge0}\frac{q^{n^2+(i-1)n}}{(q)_n}.
\]
Let $\mathcal P_1,\mathcal P_2$ be the regular languages of partitions
into parts with residues in $\Sigma_1=\{1,4\}$, respectively
$\Sigma_2=\{2,3\}$, modulo five; they have five height states, with a
$U$-loop at each allowed residue, see Figure~\ref{fig:rr-automata}.
Both pairs of automata are instances of~\S\ref{ss:compile}. The
Rogers--Ramanujan identities~\cite{Rogers1894} assert
\begin{equation}\label{eq:rr-products}
 1+\varepsilon(H_{\mathcal L_1})=\frac1{(q,q^4;q^5)_\infty},\qquad
 1+\varepsilon(H_{\mathcal L_2})=\frac1{(q^2,q^3;q^5)_\infty}.
\end{equation}
We derive~\eqref{eq:rr-products} from finite identities, by two Bailey
steps and Jacobi's triple product.

For a parameter $\gamma$, call $(\alpha,\beta)$ a Bailey pair if
\begin{equation}\label{eq:rr-bailey}
 \beta_m=\sum_{r=0}^m
 \frac{\alpha_r}{(q)_{m-r}(\gamma q;q)_{m+r}}.
\end{equation}
The unit pairs needed here have $\beta_m=\delta_{m,0}$ and
\begin{align}
 \gamma=1:&\quad \alpha_0=1,\qquad
 \alpha_r=(-1)^rq^{r(r-1)/2}(1+q^r)\quad(r>0),
 \label{eq:rr-unit-one}\\
 \gamma=q:&\quad \alpha_r=
 \frac{(-1)^rq^{r(r-1)/2}(1-q^{2r+1})}{1-q}
 \quad(r\ge0).\label{eq:rr-unit-q}
\end{align}
Both are finite identities. For \eqref{eq:rr-unit-one}, pairing the
indices $r,-r$ turns the right side of \eqref{eq:rr-bailey} into
\[
 \sum_{r=-m}^{m}
 \frac{(-1)^rq^{r(r-1)/2}}{(q)_{m-r}(q)_{m+r}}
 =\frac{(-1)^mq^{m(m+1)/2}}{(q)_{2m}}
 (q^{-m};q)_{2m}.
\]
For \eqref{eq:rr-unit-q}, pair $r,-r-1$ instead, obtaining
\[
 \sum_{r=-m-1}^{m}
 \frac{(-1)^rq^{r(r-1)/2}}{(q)_{m-r}(q)_{m+r+1}}
 =\frac{(-1)^{m+1}q^{(m+1)(m+2)/2}}{(q)_{2m+1}}
 (q^{-m-1};q)_{2m+1}.
\]
The right-hand sides vanish for $m>0$ and equal $1$ for $m=0$. These
evaluations are instances of the finite $q$-binomial identity
\[
 (y;q)_N=\sum_{j=0}^N(-1)^jq^{j(j-1)/2}
 \begin{bmatrix}N\\j\end{bmatrix}_q y^j,\qquad \begin{bmatrix}N\\j\end{bmatrix}_q=\frac{(q)_N}{(q)_j(q)_{N-j}}
\]
which follows coefficientwise from $(y;q)_{N+1}=(1-yq^N)(y;q)_N$ and the
Gaussian Pascal recurrence.

The Bailey step has a rational telescoping certificate. Put
\[
 h(m,k)=\frac{\kappa^kq^{k^2}}
 {(q)_{m-k}(q)_k(\kappa q;q)_k},\quad \zeta=q^m,\quad \xi=q^k,
\]
and
\[
 r(\zeta,\xi)=\frac{\kappa q\xi^2(1-\zeta/\xi)}{(1-q\xi)(1-\kappa q\xi)},\qquad
 \Psi(\zeta,\xi)=-\frac{\zeta(1-\xi)(1-\kappa\xi)}\xi,
\]
so that $r(\zeta,\xi)=h(m,k+1)/h(m,k)$. The rational identity
\begin{equation}\label{eq:rr-durfee-wz}
 r(\zeta,\xi)\Psi(\zeta,q\xi)-\Psi(\zeta,\xi)
 =(1-\zeta)(1-\kappa\zeta)-(1-\zeta/\xi)
\end{equation}
is checked by clearing denominators. Extend $h(m,k)$ by zero outside
$0\le k\le m$. Multiplying by $h(m,k)$ and summing
over $k$ gives $(1-q^m)(1-\kappa q^m)F_m=F_{m-1}$ for
$F_m=\sum_kh(m,k)$, because $(1-q^{m-k})h(m,k)=h(m-1,k)$; the lower flux
vanishes at $\xi=1$, the upper flux vanishes because $h(m,m+1)=0$, and
$F_0=1$. Thus
\begin{equation}\label{eq:rr-durfee}
 \sum_{k=0}^m h(m,k)=\frac1{(q)_m(\kappa q;q)_m}.
\end{equation}
The certificate $\Psi$ is found by the same search, with denominator
  $\xi$, a numerator of degree at most two in $\xi$ over
  $\Q(q,\kappa,\zeta)$, and an unknown multiple of the right side
  of~\eqref{eq:rr-durfee-wz}: clearing $\xi(1-q\xi)(1-\kappa q\xi)$ gives
  a $5\times4$ linear system with one-dimensional kernel. Substituting
\eqref{eq:rr-bailey} into
\[
 \alpha'_r=\gamma^rq^{r^2}\alpha_r,\qquad
 \beta'_m=\sum_{j=0}^m\frac{\gamma^jq^{j^2}\beta_j}{(q)_{m-j}},
\]
and applying \eqref{eq:rr-durfee} with $m\mapsto m-r$, $k=j-r$ and
$\kappa=\gamma q^{2r}$, proves that $(\alpha',\beta')$ is again a Bailey
pair.

Apply this finite transformation twice to either unit pair, and
multiply the relation at index $M$ by $(q)_M$. The result is the finite
equality
\begin{equation}\label{eq:rr-finite}
 \sum_{n=0}^M\frac{\gamma^nq^{n^2}}{(q)_n}
       \frac{(q)_M}{(q)_{M-n}}
 =\sum_{r=0}^M \gamma^{2r}q^{2r^2}\alpha_r
       \frac{(q)_M}{(q)_{M-r}(\gamma q;q)_{M+r}}.
\end{equation}
Put $\gamma=q^\delta$ with $\delta\in\{0,1\}$. Both sides converge as
$M\to\infty$, with errors of valuation at least $M+1+\delta$. On the
left, a nonconstant factorial ratio contributes at least $M-n+1$, and
\[
 n^2+\delta n+M-n+1\ge M+1+\delta\qquad(n\ge1).
\]
On the right the numerator has valuation
$w_\delta(r)=(5r^2+(4\delta-1)r)/2$ for $r>0$, and
$w_\delta(r)+M-r+1\ge M+1+\delta$. The denominator replacement at $r=0$
begins at $M+1+\delta$; the remaining denominator errors and omitted
tails have larger valuation. Hence \eqref{eq:rr-finite} gives
\begin{align*}
 \sum_{n\ge0}\frac{q^{n^2}}{(q)_n}
 &=\frac1{(q)_\infty}
   \sum_{r\in\mathbb Z}(-1)^rq^{(5r^2-r)/2},\\
 \sum_{n\ge0}\frac{q^{n^2+n}}{(q)_n}
 &=\frac1{(q)_\infty}
   \sum_{r\in\mathbb Z}(-1)^rq^{(5r^2+3r)/2}.
\end{align*}
The second equality pairs $r\ge0$ with $-r-1$, and cancels the factor
$1-q$ against $(q^2;q)_\infty$.

Jacobi's triple product follows from a finite identity:
\[
 (y;p)_N(p/y;p)_N
 =\sum_{k=-N}^N(-1)^kp^{k(k-1)/2}y^k
     \begin{bmatrix}2N\\N+k\end{bmatrix}_p.
\]
This is the preceding finite $q$-binomial theorem, after a shift of the
summation index. At $(p,y)=(q^5,q^2)$, respectively $(q^5,q^4)$, its
two sides converge with errors of valuation at least $5N+2$,
respectively $5N+1$: the Gaussian coefficient differs from
$1/(q^5;q^5)_\infty$ by a series of valuation at least $5(N-|k|+1)$, and
combining this with the quadratic theta exponent gives the stated
bounds, including the omitted tail. The resulting products are
$(q^2,q^3,q^5;q^5)_\infty$ and $(q,q^4,q^5;q^5)_\infty$. Cancellation
with $(q)_\infty$ proves \eqref{eq:rr-products}.

To bind this calculation to the quantum-plane certificates, put
$S_i=H_{\mathcal L_i}-H_{\mathcal P_i}$. Each gap language has three
interior states and each residue language five, so the normalized
difference system has eight states. Let
$\mathsf E$ be the $5\times5$ cyclic matrix with $\mathsf E_{j,j-1}=1$,
indices read modulo five; let $\I_\Sigma$ be the diagonal matrix with
entries $1_{j\in\Sigma}$, and $\chi_\Sigma$ its diagonal column. Set
\[
 \mathsf L_{\mathrm{gap}}=\begin{pmatrix}
 1-qR&0&-qR\\-q\Ups&1&0\\0&-qR&1
 \end{pmatrix},\qquad
 \mathsf L_\Sigma=\I_5-qR\mathsf E-q\Ups\I_\Sigma.
\]
The data for $S_i$ are
\[
 \mathsf L_i=\operatorname{diag}(\mathsf L_{\mathrm{gap}},\mathsf L_{\Sigma_i}),\qquad
 f_i=\binom{e_{2i-1}}{e_1},\qquad
 d_i=qR\begin{pmatrix}e_1^{\top}&-\chi_{\Sigma_i}^{\top}\end{pmatrix},
\]
where the first basis vector in $f_i$ belongs to $\Q^3$ and the
second to $\Q^5$. For the marked system use
$\mathsf L_{i,x}=\mathsf L_i|_{\Ups\mapsto x\Ups}$, $f_{i,x}=xf_i$,
$d_{i,x}=d_i$. The general construction gives $\mathsf G_b$ of size $10$
and $\mathsf G_a$ of size $18$; with their constant coordinate set to
$1$ and their state vectors the normalized solutions, their last
coordinates are the required $a,b$. Since the identities above give
$\varepsilon(S_i)=0$, Corollary~\ref{cor:indirect} yields
$\nu(a_{i,j})\ge i+j$ and $\nu(b_n)\ge n+1$, which completes both
Rogers--Ramanujan certificates.

\subsection{A global constraint: alternating sums}
For a partition $\lambda=(\lambda_1\ge\cdots\ge\lambda_\ell>0)$, put
\[
 \alpha(\lambda)=\lambda_1-\lambda_2+\lambda_3-\cdots
                 +(-1)^{\ell-1}\lambda_\ell.
\]
Bessenrodt's refinement of Euler's identity, obtained from Sylvester's
bijection \cite{Bessenrodt,CYZ}, states that for every $n,k\ge0$ the
number of distinct-part partitions of $n$ with $\alpha(\lambda)=k$ equals
the number of odd-part partitions of $n$ with exactly $k$ parts. This is
a global constraint on the partition, rather than a restriction on an
individual part or on an individual adjacent difference. Its weighted
form is
\begin{equation}\label{nl:identity}
 \sum_{\lambda\text{ distinct}}z^{\alpha(\lambda)}q^{|\lambda|}
 =\sum_{\mu\text{ odd}}z^{\ell(\mu)}q^{|\mu|}
 =\frac1{(zq;q^2)_\infty}.
\end{equation}

The refinement is still recognized by a finite weighted automaton. In
the boundary word of a distinct partition, give a letter $R$ an
additional weight $z$ precisely when an odd number of letters $U$
follow it; the total weight is $z^{\alpha(\lambda)}$. Writing
\[
 B_0=R(1-R)^{-1}U,\qquad B_1=zR(1-zR)^{-1}U,
\]
the two weighted regular series, including the empty partition, are
\begin{equation}\label{nl:languages}
 H^z_{\mathrm d}=(1+B_1)(1-B_0B_1)^{-1},\qquad
 H^z_{\mathrm o}=1+R(1-zU-R^2)^{-1}zU.
\end{equation}
Indeed, even-length distinct partitions have block pattern
$(B_0B_1)^*$, and odd-length ones have pattern $B_1(B_0B_1)^*$; the
second language allows $U$ only at odd heights, each weighted $z$.
Extracting $[z^k]$ is also a finite-state operation: record the
accumulated exponent up to $k$ and reject paths that exceed it. Thus the
identities at each fixed value of the statistic are also between
regular languages.

We prove~\eqref{nl:identity} by telescoping. For a
distinct partition of length $m$, set $\delta_i=\lambda_i-\lambda_{i+1}\ge1$,
where $\lambda_{m+1}=0$. Then
\[
 |\lambda|=\sum_{i=1}^m i\delta_i,\qquad
 \alpha(\lambda)=\sum_{\substack{1\le i\le m\\i\text{ odd}}}\delta_i.
\]
Pairing lengths $2j-1$ and $2j$ therefore gives
\[
 \frac{z^j q^{2j^2-j}}{(zq;q^2)_j(q^2;q^2)_{j-1}}
 +\frac{z^j q^{2j^2+j}}{(zq;q^2)_j(q^2;q^2)_j}
 =\frac{z^j q^{2j^2-j}}{(zq;q^2)_j(q^2;q^2)_j}.
\]
Including the empty partition and putting $\eta=zq$, $Q=q^2$, the left
side of \eqref{nl:identity} is
\[
 V(\eta)=\sum_{j\ge0}f_j(\eta),\qquad
 f_j(\eta)=\frac{\eta^j Q^{j(j-1)}}{(Q;Q)_j(\eta;Q)_j}.
\]
For $\xi=Q^j$, its two rational ratios are
\[
 r(\xi)=\frac{f_{j+1}(\eta)}{f_j(\eta)}
 =\frac{\eta\xi^2}{(1-Q\xi)(1-\eta\xi)},\qquad
 s(\xi)=\frac{f_j(Q\eta)}{f_j(\eta)}
 =\frac{\xi(1-\eta)}{1-\eta\xi}.
\]
The telescoping certificate $\Psi(\xi)=(1-\eta)(\xi-1)$ satisfies
\begin{equation}\label{nl:wz}
 r(\xi)\Psi(Q\xi)-\Psi(\xi)=(1-\eta)-s(\xi).
\end{equation}
Its lower boundary is $\Psi(1)=0$, and the upper flux at $j=N+1$ has
$\eta$-valuation at least $N+1$; after $\eta=zq$, $Q=q^2$, it has
$q$-valuation at least $2(N+1)^2-(N+1)$. Summing~\eqref{nl:wz}
therefore gives
\[
 V(Q\eta)=(1-\eta)V(\eta),\qquad V(0)=1.
\]
These equations determine $V$ uniquely in $\Q(Q)[[\eta]]$: its
coefficient of $\eta^j$ is $1/(Q;Q)_j$. The product
$(\eta;Q)_\infty^{-1}$, interpreted coefficientwise $Q$-adically,
satisfies the same equations and initial value, proving
\eqref{nl:identity}. The proof is uniform in the statistic variable $z$.

For the quantum-plane certificate, put $S=H^z_{\mathrm d}-H^z_{\mathrm o}$
and let $D(x,R)$ be its normal symbol. Set
\[
 \mathsf L_{\mathrm d}=\begin{pmatrix}
 1&0&0&-q\Ups\\
 0&1&-q\Ups&0\\
 -qR&0&1-qR&0\\
 0&-zqR&0&1-zqR
 \end{pmatrix},\qquad
 \mathsf L=\operatorname{diag}(\mathsf L_{\mathrm d},1-zq\Ups-q^2R^2),
\]
\[
 f=(0,0,q,zq,q)^{\top},\qquad
 d=(0,0,0,R,-zR).
\]
The first four states are: after a $U$ ending a $B_1$ block, after a
$U$ ending a $B_0$ block, inside a $B_0$ block, inside a $B_1$ block.
The normalized solution of $\mathsf Lm=f$ gives the marginal $C=dm$ of
$S$; its initial data at $R=0$ are
\[
 m(0)=\bigl(zq^2,q^2,q,zq,q/(1-zq)\bigr)^{\top},
\]
and all subsequent coefficients are fixed by $\mathsf L$, whose
constant-coefficient matrix at every $R$-degree is invertible over
$\Q(q,z)$. The marked system is obtained by replacing $\Ups$ by
$x\Ups$ in $\mathsf L$ and $f$ by $xf$; the output row $d$ is
unchanged. The general construction gives $\mathsf G_b$ of size $7$ and
$\mathsf G_a$ of size $12$.

The constant terms $1$ cancel in $S$, all other accepted words lie in
$R\{R,U\}^*U$, and the statistic weights lie in $\Q[z]$, so
Lemma~\ref{lem:area} and Corollary~\ref{cor:indirect} apply unchanged:
by~\eqref{nl:identity}, $\nu(a_{i,j})\ge i+j$ and $\nu(b_n)\ge n+1$.
This is a convergent holonomic certificate for every coefficient
identity in~\eqref{nl:identity}.

\subsection{Residues of parts at prescribed indices}
\label{ss:berkovich-uncu}

For a partition $\lambda_1>\cdots>\lambda_\ell>0$ into distinct parts,
let $o(\lambda)$ and $e(\lambda)$ count its odd parts at odd and even
indices, respectively. Berkovich and Uncu proved that these two
statistics have the same joint distribution as the numbers of parts
congruent to $1$ and $3$ modulo $4$, respectively
\cite[Theorem~1.1]{BerkovichUncu}. Thus
\begin{equation}\label{bu:identity}
 \sum_{\lambda\text{ distinct}}
    u^{o(\lambda)}v^{e(\lambda)}q^{|\lambda|}
 =(-q^2;q^2)_\infty(-uq;q^4)_\infty(-vq^3;q^4)_\infty.
\end{equation}
Both sides count partitions into distinct parts. The empty partition
contributes $1$ on each side and cancels in their difference. We give
the automaton kernels, and certify convergence by a recurrence for the
bounded refinement of
\cite[Theorem~2.1]{BerkovichUncu}.

\paragraph{Automata and holonomic kernels.}
The boundary word reads parts in increasing order. To recover their
indices in decreasing order, the first automaton records the parity
of the number of $U$'s still to be read. Its eight interior states are
$(h,k,r)\in\{0,1\}^3$: $h$ is the current part size modulo $2$,
$k$ is this remaining parity, and $r=1$ means that an $R$ has occurred
since the preceding $U$. Its transitions and weights are
\[
 \begin{array}{c|c|c|c}
 \text{source}&\text{letter}&\text{target}&\text{weight}\\ \hline
 (h,k,r)&R&(1-h,k,1)&1\\
 (h,k,1)&U&(h,1-k,0)&w(h,k)
 \end{array}
 \qquad
 w(h,k)=
 \begin{cases}
  1,&h=0,\\
  u,&h=1,\ k=1,\\
  v,&h=1,\ k=0.
 \end{cases}
\]
The initial $R$ enters either $(1,0,1)$ or $(1,1,1)$; a final $U$
leaves $(h,1,1)$ with weight $w(h,1)$. Exactly one initial choice
accepts any given boundary word. The flag $r$ forbids consecutive
$U$'s, hence repeated parts.

The second automaton has eight states $(h,r)\in\mathbb Z/4\mathbb Z
\times\{0,1\}$ and transitions
\[
 \begin{array}{c|c|c|c}
 \text{source}&\text{letter}&\text{target}&\text{weight}\\ \hline
 (h,r)&R&(h+1,1)&1\\
 (h,1)&U&(h,0)&\widetilde w(h)
 \end{array}
 \qquad
 \widetilde w(h)=
 \begin{cases}
 u,&h=1,\\ v,&h=3,\\1,&h=0,2.
 \end{cases}
\]
Its initial $R$ enters $(1,1)$; its final $U$ leaves any $(h,1)$
with weight $\widetilde w(h)$.

Let $\mathsf A_R,\mathsf A_U$ be the incoming transition matrices
of their disjoint union: an edge $j\to i$ contributes its weight to
entry $(i,j)$. Let $r_0$ be the column recording the three initial
entries just specified, and $d$ the row of final weights, positive for
the first automaton and negative for the second. These tables give
the marked kernel
\begin{equation}\label{bu:kernel}
 \bigl(\I-qR\mathsf A_R-qx\Ups\mathsf A_U\bigr)m_x=xr_0,
 \qquad D(x,R)=qR\,dm_x.
\end{equation}
The matrix is congruent to $\I$ modulo $q$, which selects its solution
in $\Z[u,v,x,R][[q]]$, and~\eqref{eq:canonical-boundaries} gives the
canonical holonomic $a,b$. Equivalently, the following is a kernel for
the preimage $\phi$ of Theorem~\ref{thm:A}(3d):
\begin{equation}\label{bu:phi}
 \begin{cases}
  (\I-R\mathsf A_R-q\Ups\mathsf A_U)\tilde m=r_0,\\
  (q-R)\phi=Rd\tilde m,
 \end{cases}
 \qquad b=q\Ups\phi.
\end{equation}
Here $\tilde m=m_1(R/q)$ is selected coefficientwise in $R$; the
constant matrix at degree $n$ is $\I-q^{n+1}\mathsf A_U$, with inverse
$\I+q^{n+1}\mathsf A_U$ since $\mathsf A_U^2=0$. The second equation
forces $\phi(0)=0$. The argument below proves that
$\phi\in\Z[u,v][[q,R]]$.

\paragraph{A recurrence certificate.}
Write $E_N(i,j)$ for the coefficient of $u^iv^j$ on the left of
\eqref{bu:identity} with largest part at most $2N$. Choosing which of
$2N,2N-1$ occur gives
\begin{align}\label{bu:recurrence}
 E_N(i,j)={}&E_{N-1}(i,j)+q^{2N}E_{N-1}(j,i)\notag\\
 &+q^{2N-1}E_{N-1}(j,i-1)+q^{4N-1}E_{N-1}(i,j-1).
\end{align}
The initial data are $E_0(i,j)=\delta_{i,0}\delta_{j,0}$; negative
indices give zero. Its proposed solution is
\begin{equation}\label{bu:T}
 T_N(i,j)=
 \frac{q^{2i^2-i+2j^2+j}(q^4;q^4)_N}
 {(q^4;q^4)_i(q^4;q^4)_j(q^2;q^2)_{N-i-j}}
 \quad(i,j\ge0,\ i+j\le N),
\end{equation}
and zero elsewhere. To verify it, set $X=q^{2N},Y=q^{2i},Z=q^{2j}$,
and divide each term on the right of \eqref{bu:recurrence} by
$T_N(i,j)$:
\[
\begin{array}{c|c}
\text{term}&\text{ratio}\\ \hline
T_{N-1}(i,j)&\displaystyle\frac{1-X/(YZ)}{1-X^2}\\[5pt]
q^{2N}T_{N-1}(j,i)&\displaystyle\frac{XY}{Z}\,
                                      \frac{1-X/(YZ)}{1-X^2}\\[5pt]
q^{2N-1}T_{N-1}(j,i-1)&\displaystyle\frac{X(1-Y^2)}{YZ(1-X^2)}\\[5pt]
q^{4N-1}T_{N-1}(i,j-1)&\displaystyle\frac{X^2(1-Z^2)}{Z^2(1-X^2)}
\end{array}
\]
Their sum is $1$ by the polynomial identity
\begin{equation}\label{bu:finite-check}
 (YZ-X)(Z+XY)+XZ(1-Y^2)+X^2Y(1-Z^2)=YZ^2(1-X^2).
\end{equation}
On the support boundary $i+j=N$, the factor $YZ-X$ makes the first two
terms zero; at $i=0$ or $j=0$, the factors $1-Y^2$ or $1-Z^2$ give the
required zeros. Outside the support all four terms vanish. Together
with $T_0(0,0)=1$, this proves $E_N=T_N$.

\paragraph{Convergence data.}
Put $\sigma_{i,j}=2i^2-i+2j^2+j$ and
\[
 G_{i,j}=\frac{q^{\sigma_{i,j}}(-q^2;q^2)_\infty}
                    {(q^4;q^4)_i(q^4;q^4)_j}.
\]
These are the coefficients of $u^iv^j$ in the product in
\eqref{bu:identity}: the expansion of $(-z;q^4)_\infty$ is selected
by the first-order equation $g(z)=(1+z)g(q^4z)$ and $g(0)=1$.
For $i+j\le N$, the quotient is
\[
 \frac{T_N(i,j)}{G_{i,j}}
 =\frac{(q^{2(N-i-j+1)};q^2)_\infty}
             {(q^{4N+4};q^4)_\infty}
 \equiv1\pmod{q^{2(N-i-j+1)}}.
\]
The valuation check is
\begin{equation}\label{bu:slope}
 \sigma_{i,j}-2(i+j)+1=(i-1)(2i-1)+j(2j-1)\ge0
 \qquad(i,j\in\Z_{\ge0}),
\end{equation}
each summand being a product of two integers of the same sign.
Consequently
$\nu(T_N(i,j)-G_{i,j})\ge\sigma_{i,j}+2(N-i-j+1)\ge2N+1$. When
$i+j>N$, the same inequality gives $\nu(G_{i,j})=\sigma_{i,j}\ge2N+1$.
The estimate is therefore uniform in the marking variables:
\begin{equation}\label{bu:uniform}
 \nu\!\left(\sum_{i,j}u^iv^j E_N(i,j)
                  -\sum_{i,j}u^iv^jG_{i,j}\right)\ge2N+1.
\end{equation}

By the remark after~\eqref{eq:beta}, $b_n$ is the difference between
the two weighted counts with largest part at most $n$. The second
count, truncated at $2N$, differs from its infinite product by terms
of valuation at least $2N+1$. Hence \eqref{bu:uniform} gives
$\nu(b_{2N})\ge2N+1$. Also
$\nu(b_{2N+1}-b_{2N})\ge2N+1$, since the newly added partitions
contain a part $2N+1$. Thus $b_n\to0$, and Theorem~\ref{thm:A}(3),
whose proof is unchanged for weights in $\Z[u,v]$, gives the sharper
bound and the integral preimage:
\[
 \nu(b_n)\ge n+1,\qquad
 b=q\Ups\phi,\quad \phi\in\Z[u,v][[q,R]].
\]
The same theorem supplies $\nu(a_{i,j})\ge i+j$. This completes the
convergent holonomic certificate. Its data are the two transition
tables, the rational ratios, the support and initial value, and
\eqref{bu:slope}; its main check is
\eqref{bu:finite-check}.

\paragraph{Larger moduli.}
Fu and Zeng extended~\eqref{bu:identity} as follows. Fix $L\ge2$ and
require distinct parts, with at most one part from each block
\[
 \{Lj+1,\ldots,Lj+L-1\},\qquad j\ge0.
\]
Let $o_r,e_r$ count the parts congruent to $r\pmod L$ at odd and
even indices. With
$z=\sum_{r=1}^{L-1}u_rq^r$ and $w=\sum_{r=1}^{L-1}v_rq^r$, the
one-variable case of~\cite[Theorem~2.6]{FuZeng} reads
\begin{equation}\label{bu:FuZeng}
 \sum_{\lambda}q^{|\lambda|}\prod_{r=1}^{L-1}
                 u_r^{o_r(\lambda)}v_r^{e_r(\lambda)}
 =(-q^L;q^L)_\infty(-z;q^{2L})_\infty
                         (-q^Lw;q^{2L})_\infty.
\end{equation}
For $L=2$ the block condition adds nothing to distinctness, so this
is exactly \eqref{bu:identity}.

The certificate extends to all $L$. For largest part
at most $LN$, the generating polynomial, expressed in $z,w$, satisfies
\begin{equation}\label{bu:block}
 F_N(z,w)=(1+q^{L(2N-1)}w)F_{N-1}(z,w)
       +q^{L(N-1)}(q^L+z)F_{N-1}(w,z),\quad F_0=1.
\end{equation}
There are four possibilities in the last block: neither part, its
multiple of $L$, one nonmultiple, or both. They give these four
terms, with the swap $z\leftrightarrow w$ whenever one part is added.
The coefficient of $z^iw^j$ is certified by
\begin{equation}\label{bu:TL}
 T^{(L)}_N(i,j)=
 \frac{q^{L i(i-1)+L j^2}(q^{2L};q^{2L})_N}
 {(q^{2L};q^{2L})_i(q^{2L};q^{2L})_j(q^L;q^L)_{N-i-j}},
\end{equation}
on the same support. With $X=q^{LN},Y=q^{Li},Z=q^{Lj}$,
the four ratios are exactly those in the preceding table. Thus
\eqref{bu:finite-check} verifies \eqref{bu:block} for every $L$.
The limit of \eqref{bu:TL}, summed against $z^iw^j$, is the product
in \eqref{bu:FuZeng}, by the same first-order product equations.

Since $\nu(z),\nu(w)\ge1$, each summand has valuation at least
$\sigma^{(L)}_{i,j}=Li(i-1)+Lj^2+i+j$. The convergence check becomes
\[
 \sigma^{(L)}_{i,j}-L(i+j)+L
   =L(i-1)^2+Lj(j-1)+i+j\ge1.
\]
With the tail quotient
\[
 \frac{(q^{L(N-i-j+1)};q^L)_\infty}
      {(q^{2L(N+1)};q^{2L})_\infty},
\]
this gives the uniform error $\nu(F_N-F_\infty)\ge LN+1$, and the
same argument proves convergence of the canonical boundary whenever the
marking weights have nonnegative valuation. Both sides have finite
automata: on the index side, the first automaton, with $h$ read modulo
$L$, needs one extra bit, which records whether a part not divisible by
$L$ has been used in the current block and is reset when an $R$
crosses a multiple of $L$; on the product side, each factor
$1+zq^{2Lj}$ or $1+wq^{2Lj+L}$ selects at most one part from its block.

For example, retain only residues $0,a$ at odd indices and $0,c$ at
even indices, where $1\le a,c<L$. Within the stated block condition,
these partitions are equinumerous with distinct partitions whose
parts lie in
\[
 0,\ a,\ L,\ L+c\pmod{2L}.
\]
Taking $L=4,a=1,c=3$ gives an index-dependent mod-$4$ condition
equal to the residue classes $0,1,4,7\pmod8$.

\paragraph{Arbitrary index periods.}
The automata do not depend on a product formula. For $L,M\ge1$, give
a part of residue $h\bmod L$ at index $k\bmod M$ the weight $z_{k,h}$;
weights $0$ or $1$ impose arbitrary index-dependent residue
restrictions. The states
\[
 (h,k,r)\in\Z/L\Z\times\Z/M\Z\times\{0,1\}
\]
have transitions
\[
 (h,k,r)\xrightarrow{R}(h+1,k,1),\qquad
 (h,k,1)\xrightarrow{\,z_{k,h}U\,}(h,k-1,0).
\]
The initial $R$ enters $(1,k,1)$ for every $k$, and the final $U$
leaves $(h,1,1)$ with weight $z_{1,h}$, all indices being read modulo
their periods. As before, exactly one initial choice accepts each
word, and the kernel~\eqref{bu:kernel} applies unchanged, with at most
$2LM$ interior states.

Bounded counts also have a small kernel. Let $v_N^{(k)}$ count
distinct partitions with largest part at most $N$, assigning index $k$
to the largest part and then increasing indices modulo $M$. Adding $N$
gives the recurrence
\begin{equation}\label{bu:periodic}
 v_N^{(k)}=v_{N-1}^{(k)}
            +q^N z_{k,N\bmod L}v_{N-1}^{(k+1)},
 \qquad v_0^{(k)}=1.
\end{equation}
In vector form put $C_h=\operatorname{diag}(z_{k,h})$, and let
$\mathsf S$ be the cyclic shift $(\mathsf Sv)_k=v_{k+1}$. Expand the
ordered matrix product
\[
 B(Z)=(\I+q^LZC_0\mathsf S)(\I+q^{L-1}ZC_{L-1}\mathsf S)
                   \cdots(\I+qZC_1\mathsf S)=\sum_{d=0}^L B_dZ^d,
\]
where for $L=1$ there is just the factor $\I+qZC_0\mathsf S$.
Then $f(R)=\sum_{N\ge0}v_{LN}R^N$ is selected by
\begin{equation}\label{bu:periodic-kernel}
 \left(\I-R\sum_{d=0}^L B_d\Ups^{Ld}\right)f=\boldsymbol1,
 \qquad f(0)=\boldsymbol1,
\end{equation}
since $v_{L(N+1)}=B(q^{LN})v_{LN}$. This is an $M$-component
holonomic kernel.

A matching product, and the convergence of the difference between two
proposed classes, remain separate assertions. For $L>2$ the Fu--Zeng
block condition cannot simply be dropped: without it, marking the
parts $\equiv1\pmod4$ by their index parity gives the coefficient $1+v$
at $q^3$, from $(3)$ and $(2,1)$, whereas marking the parts
$\equiv1\pmod8$ by $u$ and those $\equiv5\pmod8$ by $v$ gives $1+u$.

\subsection{An identity between five-state languages}\label{ss:nonlocal-memory}
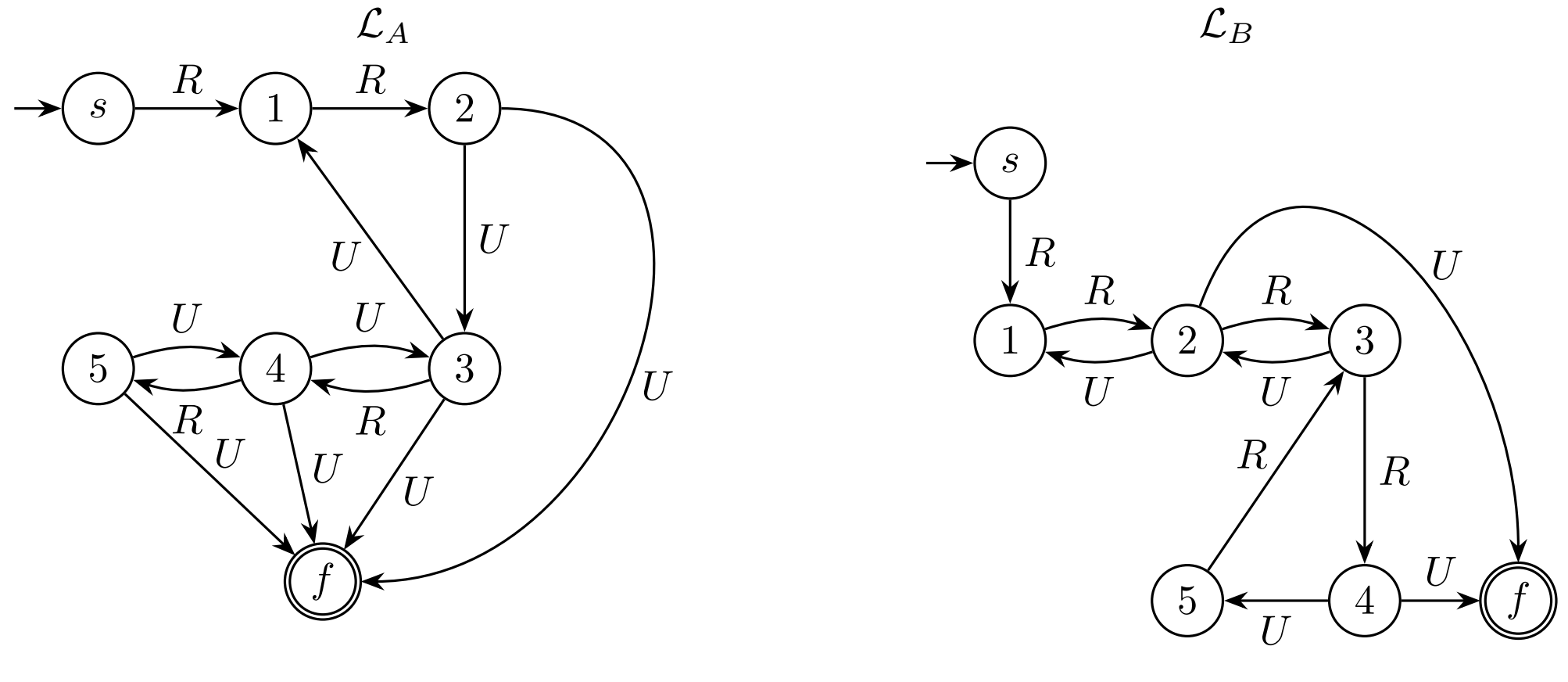
\begin{figure}[htbp]
\centering
\begingroup
\tikzset{
  cert automaton/.style={
    >=Stealth,semithick,initial text={},initial distance=4mm,
    every state/.style={draw,circle,minimum size=6mm,inner sep=1pt},
    every edge/.style={draw,->}
  },
  cert label/.style={
    postaction={decorate},
    decoration={markings,mark=at position .5 with {
      \node[anchor=center,inner sep=0pt,outer sep=0pt,
            transform shape=false] at (0,2.5mm) {$#1$};
    }}
  }
}
\begin{minipage}[t]{.49\linewidth}
\centering\vspace{0pt}
$\mathcal L_A$\par\medskip
\begin{tikzpicture}[cert automaton]
  \node[state,initial] (s) at (-1.5,2.2) {$s$};
  \node[state] (a1) at (0,2.2) {$1$};
  \node[state] (a2) at (1.6,2.2) {$2$};
  \node[state] (a3) at (1.6,0) {$3$};
  \node[state] (a4) at (0,0) {$4$};
  \node[state] (a5) at (-1.5,0) {$5$};
  \node[state,accepting] (f) at (.4,-1.8) {$f$};
  \path
    (s) edge[cert label=R] (a1)
    (a1) edge[cert label=R] (a2)
    (a2) edge[cert label=U] (a3)
         edge[out=0,in=0,looseness=1.45,cert label=U] (f)
    (a3) edge[cert label=U] (a1)
         edge[bend left=18,cert label=R] (a4)
         edge[cert label=U] (f)
    (a4) edge[bend left=18,cert label=U] (a3)
         edge[bend left=18,cert label=R] (a5)
         edge[cert label=U] (f)
    (a5) edge[bend left=18,cert label=U] (a4)
         edge[cert label=U] (f);
\end{tikzpicture}
\end{minipage}\hfill
\begin{minipage}[t]{.49\linewidth}
\centering\vspace{0pt}
$\mathcal L_B$\par\medskip
\begin{tikzpicture}[cert automaton]
  \node[state,initial] (s) at (-1.5,3.1) {$s$};
  \node[state] (b1) at (-1.5,1.6) {$1$};
  \node[state] (b2) at (0,1.6) {$2$};
  \node[state] (b3) at (1.5,1.6) {$3$};
  \node[state] (b4) at (1.5,-.6) {$4$};
  \node[state] (b5) at (0,-.6) {$5$};
  \node[state,accepting] (f) at (2.8,-.6) {$f$};
  \path
    (s) edge[cert label=R] (b1)
    (b1) edge[bend left=18,cert label=R] (b2)
    (b2) edge[bend left=18,cert label=U] (b1)
         edge[bend left=18,cert label=R] (b3)
         edge[out=70,in=90,looseness=1.6,cert label=U] (f)
    (b3) edge[bend left=18,cert label=U] (b2)
         edge[cert label=R] (b4)
    (b4) edge[cert label=U] (b5)
         edge[cert label=U] (f)
    (b5) edge[cert label=R] (b3);
\end{tikzpicture}
\end{minipage}
\endgroup
\caption{The two automata for the identity. Each has one initial
state $s$ and one final state $f$; the numbered states are interior
states.}
\label{fig:201245}
\end{figure}

Consider the two partition languages recognized by the automata of
Figure~\ref{fig:201245}. They were found by enumerating the
$14\,614\,080$ minimal deterministic interior automata with at most five
states (\S\ref{ss:searches}) and grouping them by their series to
order~$q^{50}$. Neither is defined by conditions of the
form~\eqref{eq:window}, for any $M$, $W$, $\mathcal W$ and initial
vectors, even with unbounded multiplicities. Indeed, for $n\ge W$ and
$k\equiv-2\pmod M$, the words $RRU(RU)^nU$ and $RRURRU(RU)^{n+k}URRU$ lie
in $\mathcal L_A$ but $RRU(RU)^nURRU$ does not, although each window of
the latter, with its height modulo~$M$, is a window of one of the
former; similarly, $R(RU)^nRU$ and $RR(RU)^{n+M-1}U$ lie in
$\mathcal L_B$ but $R(RU)^nU$ does not. Yet both languages count
partitions of each size equally. We give algebraic certificates of
this equality.

Write $H_A,H_B$ for their path series. Reversing a boundary
word and interchanging $R,U$ conjugates its partition, so it preserves
area. Let $H_B^\dagger$ be the series of these conjugate partitions.
We certify
\[
 S=H_A-H_B^\dagger=(1-U)a+b(1-R).
\]
The normal symbol $D(x,R)$ of $S$ has the following kernel.
Every term involving an unknown on the right has a factor $q$, so the
equations determine unique solutions in $\Q[x,R][[q]]$:
\begin{equation}\label{ex201:states}
\begin{alignedat}{3}
 m_1&=x+qx\Ups m_3,&\qquad c_1&=qx\Ups c_2,\\
 m_2&=qRm_1,& c_2&=x+qRc_1+qx\Ups c_3,\\
 m_3&=qx\Ups(m_2+m_4),& c_3&=qRc_2+qx\Ups c_4,\\
 m_4&=qRm_3+qx\Ups m_5,& c_4&=x+qRc_5,\\
 m_5&=qRm_4,& c_5&=qx\Ups c_3.
\end{alignedat}
\end{equation}
They are read directly from incoming edges: $m$ uses the first
automaton, and $c$ uses the second with all edges reversed and labels
interchanged, whose initial interior states are therefore $2$ and $4$,
and whose final interior state is $1$. The source terms $x$ and the
factors $q$ account for the two outer letters, as in
\S\ref{ss:kernels}. Thus
\begin{equation}\label{ex201:canonical}
 D(x,R)=qR(m_2+m_3+m_4+m_5-c_1),\qquad
 a(x,R)=\frac{D(x,R)-D(1,R)}{1-x}.
\end{equation}
Equations \eqref{ex201:states}--\eqref{ex201:canonical} give a
holonomic definition of $a$.

For $b$, in the operator algebra
$\Ups R=qR\Ups$, put
\begin{equation}\label{ex201:PQ}
\begin{aligned}
 P&=1-q^2(1+q)R\Ups-q^3R\Ups^2+q^8R^2\Ups^3,\\
 Q&=-1-(1+q)R-R^2+q(q^2-1)R\Ups+q^5R^2\Ups^2.
\end{aligned}
\end{equation}
Define
\begin{equation}\label{ex201:certificate}
 P\psi=q^2R,\qquad \psi(0)=0,\qquad
 \phi=Q\psi-R,\qquad b=q\Ups\phi.
\end{equation}
Since $P\equiv1\pmod R$, the first equation determines $\psi$ uniquely;
the displayed initial value is forced by it.

We now bind these data to the automata; from now on, put $x=1$
in~\eqref{ex201:states}.
Elimination gives
\begin{equation}\label{ex201:outputs}
\begin{aligned}
 qR(m_2+m_3+m_4+m_5)&=\Theta\psi,\\
 qRc_1&=(\Ups-q^5R\Ups^3)c_3-q,\\
 \Omega c_3&=q+(q-q^2)R,
\end{aligned}
\end{equation}
where
\[
 \Theta=R+qR(1-q^2R+q^2R^2)\Ups-q^5R^2\Ups^2,\qquad
 \Omega=1-q(1+q)R\Ups-q^4R\Ups^2+q^6R^2\Ups^3.
\]
For the first equality, note that the first five equations are
satisfied by
\[
 q^2R\begin{pmatrix}m_1\\m_2\\m_3\\m_4\\m_5\end{pmatrix}
 =\begin{pmatrix}
 1-q^2(1+q)R\Ups\\
 qR-q^3(1+q)R^2\Ups\\
 q^2R\Ups-q^6R^2\Ups^2\\
 q^3R^2\Ups\\
 q^4R^3\Ups
 \end{pmatrix}\psi:
\]
substituting, the first equation of~\eqref{ex201:states} becomes
$P\psi=q^2R$ and the other four hold identically. Each entry on the right side is divisible by $R$ because $\psi(0)=0$.

Put
\[
\begin{aligned}
 N&=q+q(1+q)R+q^2(1+R)\Ups-q^6R(1+R)\Ups^2,\\
 J&=q+q(1+q)R+q^2(1+q^2R)\Ups
       -q^3R(1+q^3R)\Ups^2.
\end{aligned}
\]
The two operator checks
\begin{equation}\label{ex201:intertwiner}
 \Omega N=JP,\qquad
 \Omega\bigl(q(1+R)\bigr)+J(q^2R)=q+(q-q^2)R
\end{equation}
show, by uniqueness for $\Omega\equiv1\pmod R$, that
$c_3=q(1+R)+N\psi$. Finally, direct multiplication gives
\begin{equation}\label{ex201:boundary-check}
\begin{aligned}
 &\Theta-(\Ups-q^5R\Ups^3)N-(1-R)q\Ups Q
       =\bigl(R-q^2(1+q^3R)\Ups^2\bigr)P,\\
 &q-(\Ups-q^5R\Ups^3)\bigl(q(1+R)\bigr)+(1-R)q^2R
 +\bigl(R-q^2(1+q^3R)\Ups^2\bigr)(q^2R)=0.
\end{aligned}
\end{equation}
Applying these identities to $\psi$ and using
\eqref{ex201:outputs} proves
\[
 D(1,R)=(1-R)q\Ups(Q\psi-R)=(1-R)b.
\]
Together with \eqref{ex201:canonical}, this verifies the boundary
identity. All checks are polynomial calculations using
$\Ups R=qR\Ups$; the operators $N,J$ are verification data, not
additional unknown series.

It remains to certify convergence. Every term of $1-P$ has positive
$R$-degree and coefficients in $\mathbb Z[q]$, so $P$ is an integral
kernel, and induction on the $R$-degree gives
\[
 \psi\in\mathbb Z[q][[R]],\qquad
 \phi=Q\psi-R\in\mathbb Z[q][[R]]\subset\Q[[q,R]].
\]
Together with the checks above, this is a direct certificate: $\phi$ is
the integral preimage of Theorem~\ref{thm:A}(3d), since $b=q\Ups\phi$.
Hence
\[
 \nu\bigl([R^n]b\bigr)
   =n+1+\nu\bigl([R^n]\phi\bigr)\ge n+1.
\]
Theorem~\ref{thm:A}(2) gives $\nu(a_{i,j})\ge i+j$ for the canonical
$a$ in \eqref{ex201:canonical}. Both boundaries consequently converge
$q$-adically, proving the partition identity.

\subsection{Implementation}\label{ss:implementation}
The accompanying \textsc{Julia} package~\cite{PartitionsJl}, built on
\textsc{Oscar}~\cite{OSCAR}, enumerates automata, computes their series
in~$\A$ and their marginal systems~\eqref{eq:marginal-system}, and
searches for certificates as in~\S\ref{ss:searches} and

Appendix~\ref{app:searches}. Its directory \texttt{paper} lists the
automata of the figures as transition tables, and contains one script
for each of~\S\ref{sec:teaser}
and~\S\S\ref{ss:glaisher}--\ref{ss:nonlocal-memory}. The scripts check
exactly the polynomial, rational and operator identities of the
certificates, the sizes and kernels of the searches that find them,
and, in~\S\S\ref{sec:teaser}, \ref{ss:rr} and~\ref{ss:nonlocal-memory},
that the marginal systems compiled from the automata are those displayed
at $x=1$. They also compare, to a fixed order, the series of the
automata with the stated sums and products, and the coefficients $b_n$
with the differences of bounded counts given after~\eqref{eq:beta};
these comparisons are consistency checks, on which no proof relies.
With \textsc{Julia}~1.13.1 and \textsc{Oscar}~1.8.2 the checks run in
under a minute.

\appendix
\section{Further convergence tests}\label{app:searches}
The following constructions give further sufficient tests for the
convergence of~$b$, by finite checks on holonomic data; the examples
of~\S\ref{sec:examples} do not need them.

\subsection*{A degree bound for auxiliary normalized systems}
For $\mathsf Lm=f$ with $\mathsf L\equiv\I\pmod q$ and $C=dm$ as
in~\eqref{eq:marginal-system}, choose
\[
 \kappa\ge\max(1,\deg_R\mathsf L),\qquad
 \kappa_0\ge\deg_R f+\deg_R d.
\]
Each factor $\I-\mathsf L$ raises $q$-valuation by at least one and
$R$-degree by at most $\kappa$. Thus
\[
 \deg_R[q^t]C\le\kappa_0+\kappa t,\qquad
 \nu([R^n]C)\ge
 \max\left(0,\left\lceil\frac{n-\kappa_0}{\kappa}\right\rceil\right).
\]
The same estimate applies to each component of $m$, omitting
$\deg_R d$. This supplies convergent inputs to the searches below.
For the marginal of a partition difference, Lemma~\ref{lem:area} gives
the sharper bound $\nu([R^n]C)\ge n$.

\subsection*{Polynomial and rational operator witnesses}
In~\eqref{eq:marginal-system}, seek operator rows $\Lambda_1,\Lambda_2$
with
\[
 d=(1-R)\Lambda_1+\Lambda_2\mathsf L,\qquad \Lambda_2f=(1-R)b_0,\quad b_0\in\K[R].
\]
Then $C=(1-R)\Lambda_1m+\Lambda_2f$, so $b=\Lambda_1m+b_0$. If
$\Lambda_1$ has polynomial entries with coefficients regular at $q=0$,
the preceding degree bound, applied to $m$, proves
convergence. An ansatz of bounded $\Ups$-degree for
$\Lambda_2\in\K[\Ups]^{1\times s}$ gives a finite linear system:
reduction modulo the right ideal $(1-R)\mathcal O$ is evaluation at
$R=1$ in normal order, and the remaining condition is
$(\Lambda_2f)(1)=0$. Rational or Laurent functions of $\Ups$ enlarge
the search, but an inverse shift costs valuation,
$\Ups^{-k}\sum g_nR^n=\sum q^{-kn}g_nR^n$, and poles at $\Ups=q^n$
require separately specified coefficients at finitely many exceptional
indices; such witnesses need a corresponding additional convergence
proof.

\subsection*{Contracting coefficient recurrences}
Coefficient extraction from a holonomic system gives, beyond finitely
many singular indices, a finite state recurrence $v_{n+1}=\mathsf M(q^n)v_n$
and an output $b_n=\pi(q^n)v_n$. A diagonal gauge with entries powers of
$q^n$, and a fixed block of $\ell$ steps, may make every transition
entry have valuation at least one and every output entry have valuation
bounded below. This proves linear decay along that progression by
induction. Rational functions evaluated at $q^n$ have eventually affine
valuations, by Lemma~\ref{lem:exceptional}(2); consequently
the required eventual inequalities reduce to finitely many checks. This
sufficient test does not capture cancellation special to one initial
vector, such as the one recorded by
Proposition~\ref{prop:dilation-checker}.

\subsection*{Spans of convergent holonomic series}
Call a holonomic series $g=\sum g_nR^n\in\K[[R]]$ a \emph{certified convergent input}
when it comes with selecting data and a proved bound
$\nu(g_n)\ge\delta n-c$, $\delta>0$. Normalized state systems supply
such inputs. Sums, products, sections and polynomial operator images
have effective holonomic constructions and explicit bounds. Another
useful constructor is
\[
 \nabla_\lambda g(R)=\frac{g(R)-g(R^\lambda)}{1-R},\qquad\lambda\ge2,
 \qquad\text{for which}\qquad
 [R^n]\nabla_\lambda g
 =\sum_{j=\lceil(n+1)/\lambda\rceil}^{n}g_j,
\]
so that its slope is at least $\delta/\lambda$. It is holonomic by
closure under monomial pullback and rational multiplication. For
monomial pullback, first eliminate to an equation involving only
powers of $\Ups^\lambda$, then substitute $R\mapsto R^\lambda$; this
turns $\Ups^\lambda$ into the $q$-shift of the pulled-back series.

For a finite family $g_1,\ldots,g_t$ and a degree bound, seek
\[
 C=\sum_{i,k,l}c_{i,k,l}(1-R)R^l\Ups^kg_i,\qquad c_{i,k,l}\in \K.
\]
Take
a common scalar annihilator $\Omega$ of $C$ and of the finitely many
terms on the right. By Lemma~\ref{lem:exceptional}(1), a linear
combination vanishes identically if and only if its coefficients at
$E_\Omega$ vanish, so the relation space is the kernel of the finite
matrix of those exceptional coefficients. A relation with nonzero
coefficient of $C$ expresses $b$ through certified convergent inputs, and so proves its
convergence.

Enumerating finite expressions in these constructors, with bounds on
operator degrees, gives a sound search, which terminates whenever the
canonical $b$ belongs to the span of finitely many enumerated certified convergent
inputs. In parallel, the criterion~\eqref{eq:convergence-equivalence}
eventually detects every false identity.

\section*{Acknowledgments}
I am deeply grateful to Jehanne Dousse for stimulating discussions on
partition identities, and for pointing me to many examples to which the
methods of this article apply.

ChatGPT 6 Astra and Claude Opus 5.5 helped write the code that searches
for and checks the certificates, and Claude helped prepare the
manuscript. No AI was used to set up the underlying mathematical theory.
here.

Julien David generously contributed the enumerator of finite automata
used in~\S\S\ref{ss:searches} and~\ref{ss:nonlocal-memory}
(\texttt{exhaustivautomata} in~\cite{PartitionsJl}).

\end{document}